\documentclass[twoside]{article}
\usepackage{graphicx, epstopdf, verbatim, fancyhdr, indent, wrapfig}

\usepackage{amsmath}
\input amssym.def
\input amssym.tex
\usepackage{bbm}

\skip\footins 22pt plus 2pt minus 2pt
\newdimen\jot
\def\mqth{\mathsurround=0pt }
\dimendef\dimenq=0
\def\openup{\afterassignment\qpenup\dimenq=}
\def\qpenup{\advance\lineskip\dimenq
  \advance\baselineskip\dimenq \advance\lineskiplimit\dimenq}
\def\eqalign#1{\,\vcenter{\openup1\jot \mqth
  \ialign{\strut\hfil$\displaystyle{##}$&$\displaystyle{{}##}$\hfil
  \crcr#1\crcr}}\,}
\newif\ifdtqp
\def\displqy{\global\dtqptrue \openup1\jot \mqth
  \everycr{\noalign{\ifdtqp \global\dtqpfalse
     \vskip-\lineskiplimit \vskip\normallineskiplimit
     \else \penalty\interdisplaylinepenalty \fi}}}
\def\displaylines#1{\displqy
  \halign{\hbox to\displaywidth{$\hfil\displaystyle##\hfil$}\crcr
  #1\crcr}}
\newskip\centerinq \centerinq=0pt plus 1000pt minus 1000pt

\def\bah#1{\overline#1}
\newcommand{\eps}{\varepsilon}
\newcommand{\cdt}{\mathbin{\raisebox{-.5ex}{\mbox{\huge$\cdot$}}}}
\newcommand{\ph}{\varphi}
\newcommand{\ints}{{\mathbbm{Z}}}
\newcommand{\reals}{{\mathbbm{R}}}
\newcommand{\cplxs}{{\mathbbm C}}
\newcommand{\quats}{{\mathbbm H}}
\renewcommand{\P}{{\mathbbm P}}

\newcommand{\vol}{\mathop{\rm vol}\nolimits}

\newcommand{\Lip}{\mathop{\rm Lip}\nolimits}
\newcommand{\Lk}{\mathop{\rm Link}\nolimits}
\newcommand{\IP}{\mathop{\rm IP}\nolimits}
\newcommand{\Hopf}{\mathop{\rm Hopf}\nolimits}
\newcommand{\dst}{\mathbin{{\,\displaystyle{\ast}\,}}}
\newcommand{\dast}{{\displaystyle{\ast}}}

\def\ip#1#2{\left\langle#1,#2\right\rangle}

\begin{document}
\thispagestyle{empty}
\large

\bigskip

\addtolength{\baselineskip}{4pt}

\centerline{\LARGE\bf Lipschitz minimality of Hopf fibrations}
\centerline{\LARGE\bf and Hopf vector fields\rm}

\medskip

\bigskip

\centerline{\sc Dennis DeTurck, Herman Gluck and Peter A.\ Storm\rm}

\bigskip

\large

\bigskip

\centerline{\bf Abstract\rm}

{Given a Hopf fibration of a round sphere by parallel great subspheres, we prove that the projection map 
to the base space is, up to isometries of domain and range, the unique Lipschitz constant minimizer in its 
homotopy class.

Similarly, given a Hopf fibration of a round sphere by parallel great circles, we view a unit 
vector field tangent to the fibres as a cross-section of the unit tangent bundle of the sphere, and prove
that it is, up to isometries of domain and range, the unique Lipschitz constant minimizer in its homotopy 
class.

Previous attempts to find a mathematical sense in which Hopf fibrations and Hopf vector 
fields are optimal have met with limited success.}

\medskip

\medskip

\bigskip

\large
\centerline{\bf Introduction and statement of results\rm.}

\addtolength{\baselineskip}{1pt}

The Hopf fibration $S^1\subset S^3 \to S^2$ 
of a round 3-sphere by parallel great circles was introduced by Heinz Hopf [1931]. It provided the 
first example of a homotopically nontrivial map from one sphere to another of lower dimension, 
spurring the development of both homotopy theory and fibre spaces in their infancy. Although Hopf 
first presented his map in terms of quadratic polynomials, he explained later in this paper that the 
fibres are the intersections of $S^3$ 
with the complex lines through the origin in $\reals^4=\cplxs^2$.

In his second paper on this theme, Hopf [1935] presented three families of fibrations of round 
spheres by parallel great subspheres:
$$\begin{array}{lll}
\vphantom{\displaystyle{\int}}S^1 \subset S^3\to S^2=\cplxs\P^1,& S^1\subset S^5\ \to\cplxs\P^2,\ \ldots\,,&S^1\subset
S^{2n+1}\to\cplxs \P^n,\ldots\\
S^3 \subset S^7\to S^4=\quats \P^1,& S^3\subset S^{11}\to\quats \P^2,\ \ldots\,,&S^3
\subset S^{4n+3}\to \quats \P^n,\ldots\\
\vphantom{\displaystyle{\int}}S^7 \subset S^{15}\to S^8,& & 
\end{array}$$
with base spaces the complex and quaternionic projective spaces, and with the non-associativity of 
the Cayley numbers responsible for the truncation of the third family.

This list is complete in the sense that any fibration of a round sphere by parallel great subspheres is 
isometric to one of the above (Wong [1961], Wolf [1963], Escobales [1975], Ranjan [1985]), meaning that there 
is an isometry of the total space carrying fibres to fibres.

\vfill
\eject

The isometry groups of these Hopf fibrations act transitively on the spherical total spaces, and 
so the base 
spaces inherit from them Riemannian metrics which make the projection maps into Riemannian 
submersions. In particular, the spherical base spaces $S^2$, $S^4$, and $S^8$ all have radius 
$\frac{1}{2}$.

We begin now with Theorem A.

The \it Lipschitz constant \rm $\Lip f$ of a continuous map $f\colon X\to Y$ between metric 
spaces is the smallest number $c\ge 0$ such that $d(f(x),f(x'))\le c\,d(x,x')$ for all points $x$ 
and $x'$ in $X$ . If no such
number $c$ exists, we regard the Lipschitz constant of $f$ to be infinite. Since the 
above Hopf projections are Riemannian submersions, they all have Lipschitz constant 1 .

Two maps $f_1$ and $f_2\colon X\to Y$ between metric spaces will be said to be 
\it isometric \rm if there are isometries $g_X\colon X\to X$ and $g_Y\colon Y\to Y$
such that $g_Y f_1=f_2 g_X$.

\medskip

\noindent\bf Theorem A.\ \  
\textit{\textbf{\boldmath The Lipschitz constant of any continuous map 
$$S^{2n+1}\to\cplxs \P^n\qquad\mbox{\rm or}\qquad S^{4n+3}\to \quats \P^n
\qquad\mbox{\rm or}\qquad S^{15}\to S^8$$
with nonzero Hopf invariant is $\ge 1$, and equals $1$ if and only if the map is isometric to 
the corresponding Hopf projection.}}\rm

\noindent
In particular, the Hopf projections are, up to isometries of domain and range, the unique 
Lipschitz constant 
minimizers in their homotopy classes.

The proof is entirely elementary metric geometry. Aiming for clarity of presentation, we give the argument 
first in the case of $S^1\subset S^3\to S^2$ , then in the case of 
$S^1\subset S^5\to\cplxs\P^2$, 
and after that explain the minor adjustments needed to carry out the proof in general.

\medskip

We turn to Theorem B.

Let $S^n$ denote the round n-sphere of radius 1 , let $S^n\times S^n$ be given the 
product metric, let $\Delta S^n=\{(x,x)\,:\,x\in S^n\}\subset S^n\times S^n$ be the diagonal, 
which is isometric to a round $n$-sphere of radius $\sqrt{2}$, 
and let $i\colon \Delta S^n\to S^n\times S^n$ denote the inclusion map.

\medskip

\noindent\bf Theorem B. \ \ 
\textit{\textbf{\boldmath The Lipschitz constant of any map 
$\Delta S^n\to S^n\times S^n$ which is homotopic to the inclusion is $\ge 1$, 
and equals $1$ if and only if the map is isometric to the inclusion.}}\rm

\noindent
This result can be appreciated by contrasting it with the following facts, pointed out to us some 
time ago by Walter Wei [1985].
\begin{enumerate}
\item[(1)] The diagonal circle $i\colon\Delta S^1\to S^1\times S^1$ is 
length-minimizing in its homology class, and any other length-minimizer in that class is isometric to it.
\item[(2)] The diagonal 2-sphere $i\colon \Delta S^2\to S^2\times S^2$ 
is area-minimizing in its homology class, but there are other area-minimizers, such as 
$S^2\vee S^2$, in the same class.
\item[(3)] For $n\ge 3$, the diagonal $n$-sphere $i\colon \Delta S^n\to S^n\times S^n$ is 
not volume-minimizing in its homology class, since $S^n\vee S^n$ 
lies in the same class but has smaller volume.
\end{enumerate}
Thus, minimizing the ``stretch'' (Lipschitz constant) of a map in its 
homotopy class may be viewed as an alternative to minimizing the area or volume of a cycle in its 
homology class, and yields different results.

We turn to Theorem C.

Given a Hopf fibration of $S^{2n+1}$ by parallel great circles, let $v$ denote either 
of the two unit vector fields on $S^{2n+1}$ which are tangent to these fibres. Then 
define $V\colon S^{2n+1}\to US^{2n+1}$ by $V(x) = (x, v(x))$, so that $V$ is the 
corresponding cross-section of the unit tangent bundle of $S^{2n+1}$. We will also 
refer to the image $V(S^{2n+1})$ as a ``Hopf vector field'', and let 
$i\colon V(S^{2n+1})\to US^{2n+1}$ denote the inclusion map.

\medskip

\noindent \bf Theorem C. \ \  
\textit{\textbf{\boldmath The Lipschitz constant of any map of the Hopf vector 
field $V(S^{2n+1})$ into the unit tangent bundle $US^{2n+1}$ which is homotopic to the 
inclusion is $\ge 1$, and equals $1$ if and only if the map is isometric to the inclusion.}}\rm

\noindent
In other words, Hopf vector fields are, up to isometries of domain and range, the unique Lipschitz 
constant minimizers in their homotopy classes.
This theorem, which asserts the Lipschitz minimality of Hopf vector fields on spheres, may be compared
with attempts to prove their volume-minimality within the unit tangent bundle, as follows.
\begin{enumerate}
\item On the 3-sphere, the 3-cycle $V(S^3)\subset US^3$	is volume-minimizing in its 
homology class. This was shown by a calibrated geometry argument in Gluck and Ziller [1986] .
\item On the 5-sphere, the 5-cycle $V(S^5)\subset US^5$	 is not volume-minimizing 
in its homology class, and indeed, is not even a local minimum. This was shown 
by David Johnson [1988], and likewise on all higher odd-dimensional spheres.
\item Sharon Pedersen [1993] showed that on each odd-dimensional sphere, beginning with $S^5$, 
there exist unit vector fields of exceptionally small volume which converge to a vector field with 
one singularity. She conjectured that on these spheres there are no unit vector fields of minimum 
volume at all, but that instead her limiting vector-field-with-singularity has minimum volume in its 
homology class in the unit tangent bundle. To support this, she showed that, as the ambient 
dimension increases, the volumes of her singular fields grow at the same rate as the known lower bound
 for volumes of nonsingular vector fields.
\end{enumerate}
\textit{\textbf{Two remarks.}}
\begin{itemize}
\item If the smooth submanifold $M$ of the Riemannian manifold $N$ is a volume-minimizing cycle in 
its homology class, then the inclusion map $i\colon M\to N$ is a Lipschitz constant minimizer in its homotopy class.
\item Theorem C will follow quickly from Theorem B.
\end{itemize}

We conclude with Theorem D.

We suspect that many natural geometric maps, such as Riemannian submersions of compact
homogeneous spaces, are Lipschitz constant minimizers in their homotopy classes, unique up 
to isometries of domain and range.

We give one further example of this in the theorem below.

Let $V_2\reals^4$ be the Stiefel manifold of orthonormal 2-frames in 4-space, with the metric 
inherited from its natural inclusion in $S^3\times S^3$, and let $G_2\reals^4$ be the Grassmann 
manifold of oriented 2-planes through the origin in 4-space. The natural projection map 
$V_2\reals^4	\to G_2\reals^4$ takes an orthonormal 2-frame to the 2-plane oriented by this 
ordered basis, and has Lipschitz constant 1 with respect to the Riemannian submersion metric that 
it induces on the Grassmann manifold.

\medskip

\noindent\bf Theorem D. \ \  
\textit{\textbf{\boldmath The Lipschitz constant of any map of 
$V_2\reals^4\to G_2\reals^4$ homotopic to the Stiefel projection is $\ge 1$, with equality if and only 
if the map is isometric to this projection.}}\rm

\noindent
To prove this theorem, we will observe within the Stiefel projection $V_2\reals^4\to G_2\reals^4$ 
two families of Hopf projections $S^3\to S^2$, whose Lipschitz minimality, unique up to isometries 
of domain and range, was established in Theorem A. They provide the framework for the proof.

\vfill
\eject

\noindent \bf Acknowledgements.\rm

We are especially grateful to Shmuel Weinberger, whose interest in determining the minimum 
Lipschitz constant $L(d)$ for maps from $S^3\to S^2$ of Hopf invariant $d$, and the 
asymptotic behavior of this function for large $d$, got us into this subject in the first place.

Many thanks to Dennis Sullivan for the proof of the Key Lemma used here in the argument for 
Theorem B, and also to Rob Schneiderman and David Yetter for the main ideas in an earlier proof of 
this lemma, centered on immersions of $n$-manifolds in $2n$-space.

Special thanks to Olga Gil-Medrano, Paul Melvin, Shea Vela-Vick and Clayton Shonkwiler, as well as 
to Kerstin Baer, Eric Korman and Haomin Wen, for their substantial help during the preparation of this paper.

Storm was partially supported by NSF grant DMS-0904355 and the Roberta and Stanley Bogen 
Visiting Professorship at Hebrew University.

\bigskip

\noindent \bf Harmonic maps.\rm

We noted above that, beginning on $S^5$, Hopf vector fields are no 
longer volume-minimizing cycles in their homology classes in the unit tangent bundle. So it is natural 
to ask if they might be energy-minimizers there.

If $L\colon V\to W$ is a linear map between inner product spaces, its \it energy \rm 
$\Vert L\Vert^2$ is defined to be the sum of the squares of the entries in a matrix for $L$ 
with respect to orthonormal bases for both $V$ and $W$, and is easily 
checked to be independent of such choices.

The \it energy \rm of a smooth map $f \colon M\to N$ between Riemannian manifolds (with $M$ compact) 
is then defined by
$$E(f ) = \frac{1}{2}\int_{x\in M} \Vert df_x\Vert^2\,d(\vol).$$  
Such a map is said to be \it harmonic \rm if it is a critical point of the energy function,
that is, if
$$ \left.\frac{dE(f_t)}{dt}\right\vert_{t=0} = \left.\frac{d}{dt}\right\vert_{t=0}
\frac{1}{2}\int_{x\in M} \Vert d(f_t)_x\Vert^2\,d(\vol)= 0$$ 
for all one-parameter families $\{f_t\}$ of maps from $M\to N$ with $f_0 = f$.

\vfill
\eject

Hopf projections are harmonic maps (Fuller [1954], Eells and Lemaire [1978]); unfortunately, harmonic 
maps from spheres to compact Riemannian manifolds are always unstable (Xin [1980]).

If a vector field $V$ on a Riemannian manifold $M$ is regarded as a map of $M$ to its 
tangent bundle $TM$, 
then $V$ is harmonic if and only if it is parallel (Nouhaud [1977], Ishihara [1979], Konderak [1992]).

By contrast, if a unit vector field on $M$ is regarded as a map into its unit tangent bundle $UM$ 
with the standard Sasaki metric, then Hopf vector fields $V_H$ on all odd-dimensional spheres 
are unstable harmonic maps. On $S^3$	 there are no other unit vector fields which are 
harmonic (Han and Yim [1996]).

If we now only look at cross-sections of the unit tangent bundle $UM$, rather than at all 
maps of $M\to UM$, then the Hopf vector fields $V_H\colon S^n\to US^n$ are still 
unstable for $n = 5, 7, 9,\ldots$ (Wood [1997]). But for $n = 3$ they are stable, and in fact 
local minima of the energy (Wood [1999]).

The relation between volume and energy of unit vector fields on spheres and related spaces has been studied over the past decade by Olga Gil-Medrano and her collaborators. A cross-section of their papers is listed in the references.

\vfill
\eject

\centerline{\bf PART I. PROOF OF THEOREM A FOR MAPS FROM $S^3$ TO $S^2$.}

\noindent\bf Linking.\rm

Since the Hopf invariant of a map reports linking of inverse images, we begin by commenting on this 
from two perspectives, homology and cohomology.

\noindent \textit{\textbf{Homology}}\rm.\ \  Let $K$ and $K'$ be disjoint oriented smooth simple closed 
curves in $\reals^3$. Let $S$ and $S'$ be oriented surfaces bounded by $K$ and $K'$, in 
general position with respect to one another. Then the linking number $\Lk(K, K')$ of $K$ and 
$K'$ can be defined to be the oriented intersection number of $K$ with $S'$ or of $K'$ with $S$, 
and standard arguments show that both quantities are equal, and hence independent of the choices 
of $S$ and $S'$.

\noindent\textit{\textbf{Cohomology}}\rm.\ \  Given $K$ and $K'$ as above, they have disjoint open 
tubular neighborhoods $U$ and $U'$, each an open solid torus. By Poincar\'e duality, the 
one-dimensional homology of $U$ is isomorphic to its two-dimensional cohomology with compact 
support, $H_1(U; \ints) \cong H^2_c(U; \ints)$, and likewise for $U'$. 
Let $\beta$ and $\beta'$ be 2-forms with compact support in $U$ and $U'$ which are 
dual in this way to $K$ and $K'$.

Extend $\beta$ and $\beta'$ over $\reals^3$ to be zero outside $U$ and $U'$, and then let 
$\alpha$ and $\alpha'$ be 1-forms with compact support in $\reals^3$ such that 
$d\alpha = \beta$ and $d\alpha'=\beta'$. Then we can define
$$\Lk(K,K') = \int_{\reals^3}\alpha\wedge\beta' = \int_{\reals^3}\alpha'\wedge\beta,$$
and standard arguments show that both integrals are equal, hence independent of the 
choices of $\alpha$ and $\alpha'$, and that this definition of linking number coincides with the one 
given above.

\medskip

\noindent\bf The Hopf invariant of a map from $S^3$ to $S^2$.\rm

We give two equivalent definitions of the Hopf invariant of a continuous map 
$f\colon S^3\to S^2$, and refer the reader to Bott and Tu [1982, pp. 227-239] for further details.

\noindent
\bf(1)\rm\ \  Homotope $f$ to a smooth map, which we still call $f$, and take any 
two regular values $y$ and $y'$. Then the inverse images $f^{ -1}(y)$ and $f^{-1}(y')$ 
are smooth 1-dimensional submanifolds of $S^3$, hence each is a finite union, say $K$ and $K'$, 
of smooth simple closed curves, which we orient as follows. Start with orientations of the domain $S^3$
and the range $S^2$. Suppose $x$ is a point of $K = f^{-1}(y)$. Choose a small disk 
$D_x$ in $S^3$ through $x$, transverse there to $K$ . Orient $D_x$ so that the restriction of 
$f$ to it is orientation-preserving. Then orient the component $K_x$ of $K$ containing $x$ so 
that the orientation of $D_x$ followed by the orientation of $K_x$ agrees with the orientation of $S^3$.
Continue in this way to orient all the components of $K$ and $K'$. Then define the Hopf invariant of 
$f$ to be the total linking number of all the components of $K$ with all the components of $K'$.

Hopf [1931] showed that this definition is independent of the choice of regular values $y$ and $y'$ 
of $f$, and that it depends only on the homotopy class of $f$, not on the particular choice of $f$ itself.

\noindent
\bf(2)\rm\ \  Use singular cohomology with integer coefficients, and let $\omega$ be a 
2-dimensional cocycle on $S^2$ with $\ip{\omega}{S^2} = 1$. Then the pullback 
$f^\dast \omega$ is a 2-cocycle on $S^3$. Since $H^2(S^3; \ints) = 0$, there is a 
1-dimensional integral cochain $\alpha$ on $S^3$ such that $d\alpha = f^\dast\alpha$. 
Then the integer $\ip{\alpha\cup f^\dast\omega}{S^3}$
is defined to be the Hopf invariant of $f$. Note that we are using ``\,$d$\,'' instead of ``\,$\delta$\,''  for the coboundary map, as in the case of differential forms.

One shows that this definition is independent of the choice of 2-cocyle $\omega$ on $S^2$, and of 
the choice of 1-cochain $\alpha$ on $S^3$ such that $d\alpha = f^\dast \omega$, and that it depends 
only on the homotopy class of $f$, not on the particular choice of $f$ itself.

Unlike Hopf 's definition, this one does not require us to first homotope $f$ to make it smooth. 
However, if $f$ is smooth, we can use de Rham cohomology, let $\omega$ be a smooth 2-form 
on $S^2$ such that $\int_{S^2}\omega=1$, let $\alpha$ be a smooth 1-form on $S^3$ such that 
$d\alpha= f^\dast\omega$, and then the integral $\int_{S^3}\alpha\wedge f^\dast\omega$
gives the Hopf invariant of $f$. This is the approach of J.H.C.Whithead [1947], who showed it 
to be equivalent to Hopf's.

\medskip

\noindent\bf Mix-and-match formula for the Hopf invariant.\rm

Since we will be looking at all continuous maps $f\colon S^3\to S^2$, not known in advance to be 
smooth, we favor Whitehead's approach to the Hopf invariant, phrased as above in the language 
of singular cohomology with integer coefficients.

Here is a curiosity of that approach. Initially it is just a play with two actors: the 2-dimensional 
cocycle $\omega$ on $S^2$ with $\ip{\omega}{S^2}=1$,
and the 1-dimensional cochain $\alpha$ on $S^3$ with $d\alpha = f^\dast\omega$, with the 
Hopf invariant of $f$ given by
$$\Hopf(f) = \ip{\alpha\cup f^\dast\omega}{S^3}.$$
A third actor can be introduced: another 2-dimensional cocycle $\omega'$ on $S^2$
with $\ip{\omega'}{S^2}=1$,
and then we claim that 
$$\Hopf(f) = \ip{\alpha\cup f^\dast\omega'}{S^3},$$
which we view as a ``mix-and-match'' formula. To verify its correctness, note that $\omega$ 
and $\omega'$ are cohomologous on $S^2$, so we can write $\omega-\omega'=d\eta$,
for some 1-cochain $\eta$ on $S^2$. Then
$$\alpha\cup f^\dast\omega-\alpha\cup f^\dast\omega'=\alpha\cup f^\dast d\eta=\alpha\cup df^\dast\eta.$$
Now
$$d(\alpha\cup d^\dast\eta)=d\alpha\cup f^\dast\eta-\alpha\cup df^\dast\eta.$$
So integration by parts yields
$$\eqalign{
\ip{\alpha\cup f^\dast\omega}{S^3}-\ip{\alpha\cup f^\dast\omega'}{S^3}&=
\ip{\alpha\cup df^\dast\eta}{S^3}\cr
&=\ip{d\alpha\cup f^\dast\eta}{S^3}-\ip{d(\alpha\cup d^\dast\eta)}{S^3}\cr
&=\ip{d\alpha\cup f^\dast\eta}{S^3},
}$$
since $\ip{d(\alpha\cup f^\dast\eta)}{S^3}=0$ by Stokes's theorem.

But $d\alpha = f^\dast\omega$, and hence
$$d\alpha\cup f^\dast\eta=f^\dast\omega\cup f^\dast\eta=f^\dast(\omega\cup\eta)=0,$$
since $\omega\cup\eta$ is a 3-form on $S^2$, and hence identically zero. 
This verifies the mix-and-match formula above.

\medskip

\noindent\bf A sufficient condition for the Hopf invariant to be zero.\rm

We put the mix-and-match formula to immediate good use.

As motivation, suppose that $f\colon S^3\to S^2$ is a smooth map, with $y$ and $y'$ as regular 
values, so that the Hopf invariant of $f$ is given by the formula
$$\Hopf(f) = \Lk(K,K')$$ 
where $K = f^{-1}(y)$ and $K' = f^{-1}(y')$ are smooth oriented links in $S^3$.

Suppose there is an open set $U$ in $S^3$ which contains $K$, excludes $K'$, and has trivial 
1-dimensional homology:
$$K \subset U \subset S^3 - K'	\quad \mbox{\rm and}\quad H_1(U; \ints) = 0.$$ 
Then the link $K$ bounds a 2-chain $S$ in $U$, automatically disjoint from $K'$,
and hence the linking number of $K$ and $K'$ must be zero. Thus $\Hopf(f ) = 0$.

The following version of this, which applies to continuous rather than smooth maps, is suitable for our purposes.

\vfill
\eject

\noindent\bf Lemma E (Preliminary version)\it.\ \  Let $f\colon S^3\to S^2$ be a continuous map, 
and let $y$ and $y'$ be two points of $S^2$, with inverse images 
$K = f^{-1}(y)$ and $K' = f^{-1}(y')$. Suppose there is an open set $U$ in $S^3$ such that
$$K\subset  U \subset \bah{U} \subset S^3-K' \quad \mbox{\rm and}\quad
H_1(U; \ints) = 0.$$ 
Then the Hopf invariant of $f$ is zero.\rm

\medskip

\noindent\bf Comment\rm.\ \  In the above statement, the symbol $\bah{U}$ denotes the closure 
of $U$, and if the chain of inclusions holds, we will say that $U$ \it separates \rm $K$ from $K'$. 
Note that if $U$ separates $K$ from $K'$, then $S^3-U$ separates $K'$ from $K$ .

\medskip

\noindent\it Proof\rm.\ \  First we need to refine the above chain of inclusions by finding small 
open sets $V$ and $V'$ about $y$ and $y'$ in $S^2$ so that
$$f^{-1}(V)\subset U \subset S^3 - f^{-1}(V').$$
To find $V$, note that the image under $f$ of the compact set $S^3-U$ is compact and 
hence closed in $S^2$, and misses the point $y$ because $f^{-1}(y)\subset U$. Therefore 
$V = S^2-f(S^3-U)$ is an open neighborhood of $y$ in $S^2$ whose inverse image 
$f^{-1}(V)$ lies in $U$, as desired. To find $V'$, repeat this with $S^3-\bah{U}$ in place of 
$U$ and $y'$ in place of $y$.

\begin{figure}[h!]
\center{\includegraphics[height=140pt]{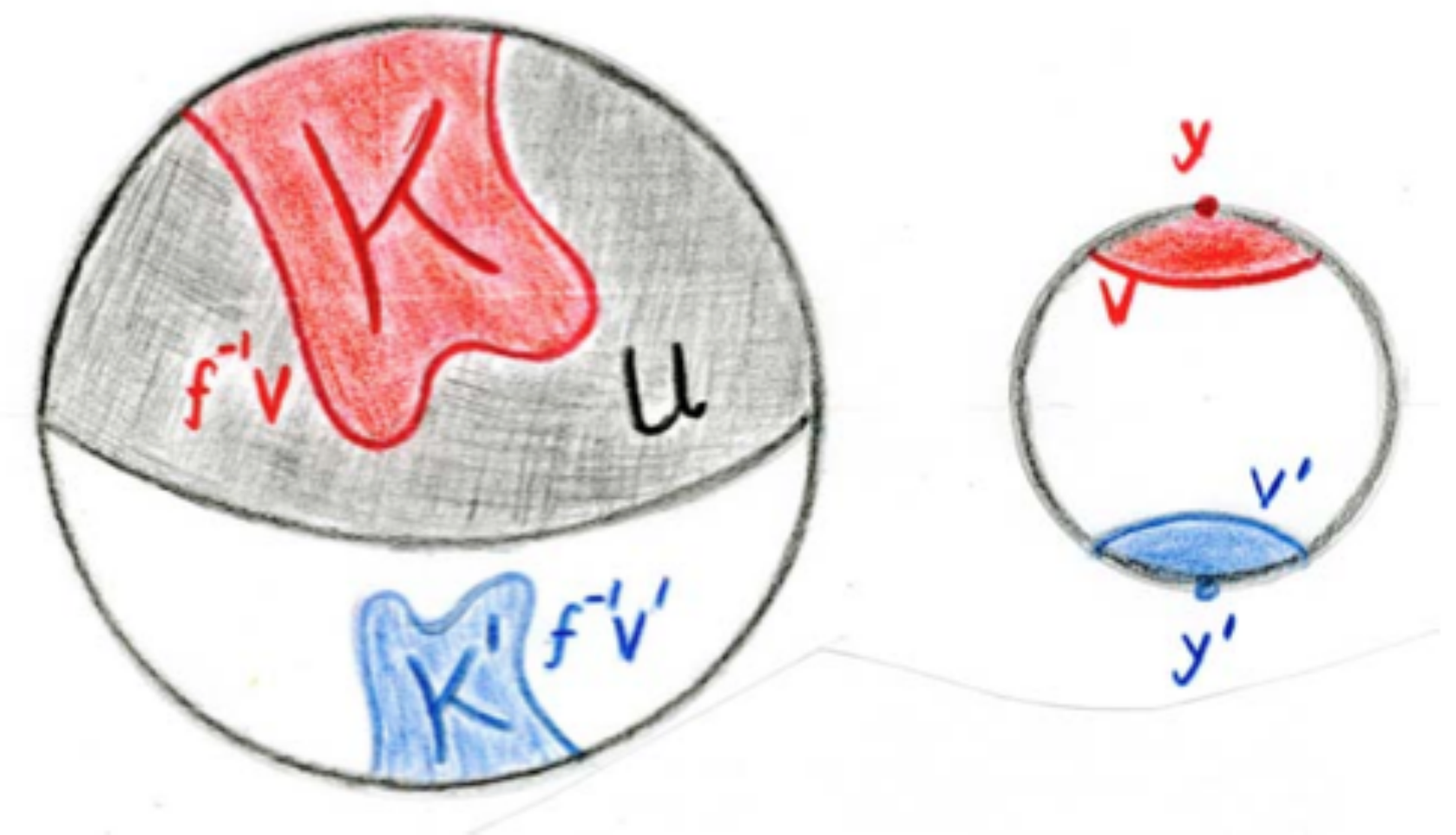}
\caption{\large{\textit{\textbf{\boldmath $f\colon S^3\to S^2$}}}}
}
\end{figure}

Now let $\omega$ be a 2-dimensional singular cocycle on $S^2$ with support
in $V$, such that $\ip{\omega}{S^2}=1$, and likewise for $\omega'$ and $V'$.

Then $f^\dast\omega$ is a 2-dimensional singular cocycle on $S^3$ with support in \linebreak
$f^{-1}(V)\subset U$. By Poincar\'e duality, the 2-dimensional singular cohomology of $U$ with 
compact supports, $H^2_c(U; \ints)$, is isomorphic to $H_1(U; \ints)$, which by hypothesis is 
zero. Hence there is a 1-dimensional cochain $\alpha$ on $S^3$ with compact support inside $U$, 
such that $d\alpha = f^\dast\omega$.

Now $\alpha$ and $f^\dast\omega'$ have supports inside the disjoint open sets $U$ and 
$f^{-1}(V')$, and therefore the cohomology class $[\alpha\cup f^\dast\omega'] = 0$. 
Then by the mix-and- match formula for the Hopf invariant, we have
$$\Hopf(f) = \ip{\alpha\cup f^\dast\omega'}{S^3}= 0.$$

\bigskip

\noindent\bf Plan of the proof of Theorem A for maps of $S^3$ to $S^2(1/2)$.\rm

We will show that any continuous map $f\colon S^3\to S^2(1/2)$ with nonzero Hopf invariant 
has Lipschitz constant $\ge 1$, with equality if and only if the map is isometric to the Hopf projection.
There are four steps to the proof, as follows.

\noindent \bf Step 1\rm.\ \  We show that for each point $y\in S^2$, its inverse image 
$f^{-1}(y)$ lies on some great 2-sphere in $S^3$.

\noindent\bf Step 2\rm.\ \ We show that each inverse image $f^{-1}(y)$ is a great circle in $S^3$ . 

\noindent\bf Step 3\rm.\ \ We show that any two such great circles $f^{-1}(y)$ and 
$f^{-1}(y')$ are
parallel to one another.
 
\noindent\bf Step 4\rm.\ \  We conclude that $f$ is isometric to the Hopf projection.

In what follows, we use the phrases ``fibre of $f$'' and ``point-inverse-image of $f$'' interchangeably.

\medskip

\noindent\bf Step 1. Each fibre of $f$ lies on a great 2-sphere in $S^3$.\rm

Let $f\colon S^3 \to S^2(1/2)$ be a map with nonzero Hopf invariant and with 
Lipschitz constant $\le 1$.

If $A$ is a subset of $S^3$ and $r$ is a positive real number, $N(A, r)$ will denote the open 
$r$-neighborhood of $A$,
$$N(A,r) = \{p\in S^3\,:\, d(p,A) < r\}.$$
We begin the argument by choosing at random a point $y\in S^2(1/2)$, and letting 
$K = f^{-1}(y)$ denote its inverse image in $S^3$. Since $\Lip f\le1$, no point in 
$N(K, \pi/2)$ can map to the antipodal point $-y$ in $S^2(1/2)$.

On the other hand, some point in $S^3$ must map to $-y$ because $f$ is homotopically nontrivial, 
and hence onto. Say $f(-x) = -y$.

Since $\Lip f\le 1$, the point $-x$ can not lie in $N(K,\pi/2)$, and therefore no point of $K$ can lie 
in the open hemisphere $N(-x,\pi/2)$. Hence $K$ must lie in the closed hemisphere of 
$S^3$ centered at $x$, as shown below.

We depict $y$ and $-y$ as north and south poles of $S^2(1/2)$, and $x$ and $-x$ as north and 
south poles of $S^3$, with $ES$ as the corresponding equatorial great 2-sphere.

\begin{figure}[h!]
\center{\includegraphics[height=180pt]{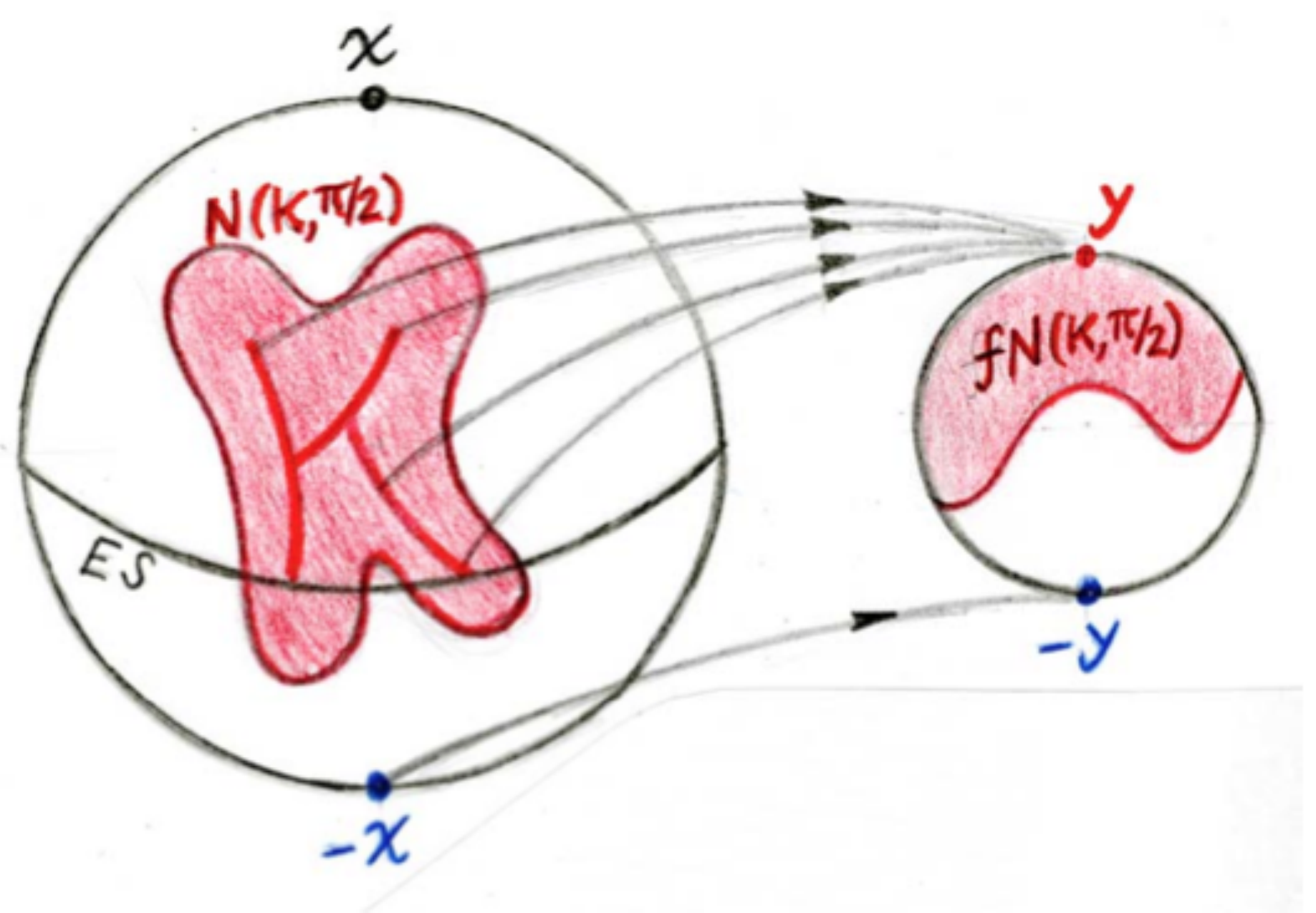}
\caption{\large{\textit{\textbf{\boldmath $f\colon S^3\to S^2(1/2)$}}}}
}
\end{figure}

The above figure shows $x$ lying outside $N(K, \pi/2)$, and we argue now that this is correct.

Suppose to the contrary that $x$ lies inside $N(K, \pi/2)$. Since $K$ lies in the closed 
northern hemisphere of $S^3$ centered at $x$, we know that the half open geodesic arc 
$[p, x)$ from each point $p$ of $K$, up to but not including $x$, must lie in $N(K, \pi/2)$. If $x$ also 
lies in $N(K, \pi/2)$, then each closed geodesic arc $[p, x]$ lies in $N(K, \pi/2)$.

Thus $K$ can be contracted along these geodesic arcs within $N(K, \pi/2)$ to the single point $x$. 
If $f$ were smooth with regular values at $y$ and $-y$, this would be enough to show that the 
linking number of $K = f^{-1}(y)$ and $K' = f^{-1}(-y)$ is zero, and hence that the Hopf invariant 
of $f$ is zero. This contradiction would then show that $x$ must indeed lie outside 
$N(K, \pi/2)$, confirming the accuracy of the above figure.

But we don't know in advance that $f$ is smooth, and so must work a little harder to expose the contradiction.

Consider our assumption (contrary to fact) that $x$ lies in the open set $N(K, \pi/2)$. Then for 
some small $\eps> 0$, the closure of the $3\eps$-ball $N(x, 3\eps)$ also lies in 
$N(K, \pi/2)$. It follows that

\noindent
(1)\ \  The closure of the $2\eps$-ball $N(x, 2\eps)$ lies in $N(K, \pi/2-\eps)$.
 
Letting $K' = f^{-1}(-y)$, and noting our assumption that $\Lip f\le 1$, we have 

\noindent
(2)\ \  $N(K, \pi/2-\eps)$ and $N(K',\eps)$ must be disjoint.

Now let $C$ denote the cone over the open set $N(K,\eps)$ from the north pole $x$ of $S$, 
meaning the union of all geodesic arcs from points of $N(K,\eps)$ to $x$. 
Denote such a geodesic arc by $[p, x]$, and note that it has length less than $\pi/2+\eps$. If we 
stop that geodesic arc $2\eps$ short of $x$, say at the point $x'$, then the subarc $[p, x']$ lies 
entirely in $N(K, \pi/2-\eps)$.

We can complete the trip along the geodesic arc from $x'$ to $x$ within the closure of the ball 
$N(x, 2\eps)$ , and hence by (1) above within the open set $N(K, \pi/2-\eps)$.

Now let
$$U = C\cup N(x,2\eps),$$
the union of two cones in $S^3$ with vertices at $x$ . Since $C$ is a cone over the open set 
$N(K,\eps)$, it is open at all of its points, save possibly at $x$. Addition of the open set 
$N(x, 2\eps)$ repairs this deficit, and so the set $U$ is open. As the union of two cones, 
it is contractible within itself to $x$.

\begin{figure}[h!]
\center{\includegraphics[height=140pt]{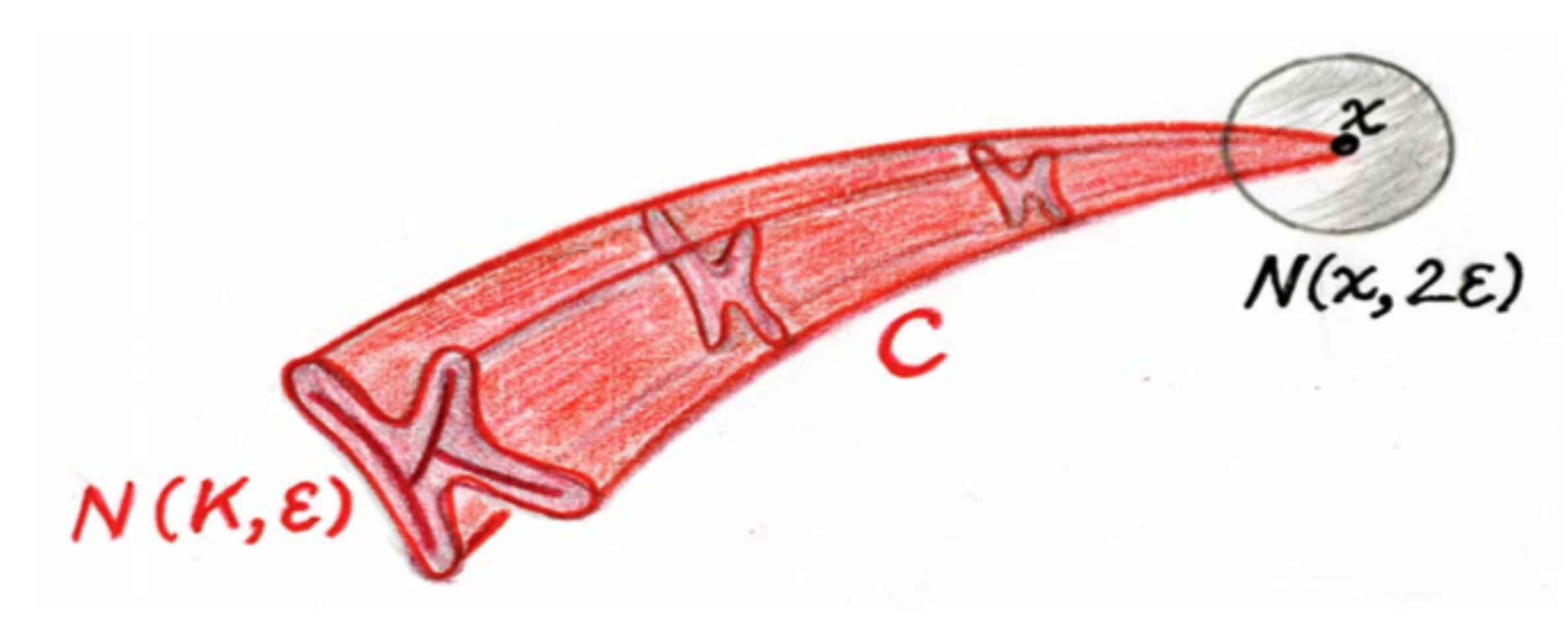}
\caption{\large{\textit{\textbf{\boldmath The contractible open set $U=C\cup N(x,2\eps)$}}}}}
\end{figure}

By construction, we have
$$K\subset U\subset N(K,\pi/2-\eps),$$ 
which is disjoint from $N(K', \eps)$ by (2) above. Hence
$$K\subset U\subset\bah{U}\subset S^3- K'.$$ 
Since $U$ is contractible, it certainly satisfies $H_1(U; \ints) = 0$.

Thus the conditions of Lemma E are satisfied, and we conclude that the Hopf invariant of $f$ is zero.

This contradiction shows that the north pole $x$ can not lie inside $N(K,\pi/2)$, and confirms the 
accuracy of its placement in Figure 2.

Since $x$ cannot lie inside $N(K, \pi/2)$, it follows that no point of $K$ can lie in the open northern
hemisphere $N(x, \pi/2)$. Since we already know that $K$ lies in the closed northern hemisphere, 
it follows that $K$ must lie on its boundary, the equatorial\linebreak
 2-sphere $ES$, completing Step 1.

\medskip

\noindent\bf Step 2. Each fibre of $f$ is a great circle in $S^3$.\rm

So far, we know that the fibre $K = f^{-1}(y)$ lies on the equatorial 2-sphere $ES$, and we intend 
to recreate there the same situation we had on the full 3-sphere $S^3$.

To begin, some point of the fibre $K' = f^{-1}(-y)$ must also lie on $ES$. Otherwise, for a sufficiently 
small positive value of $\eps$, the open equatorial region $U = N(ES,\eps)$ would separate $K$ from $K'$. 
Since $H_1(U;\ints) = 0$, Lemma E would imply that $\Hopf(f ) = 0$.

Stealing notation from the previous section, let $-x$ now denote a point of $ES$ with $f(-x) = -y$, 
so that $-x$ lies in the fibre $K' = f^{-1}(-y)$. Note that this point $-x$ is entirely different from 
the point of the same name in the previous section.

Since $\Lip f \le 1$, the point $-x$ can not lie in $N(K, \pi/2)$, and therefore no point of $K$ can lie in 
the open 2-dimensional hemisphere $ES \cap N(-x, \pi/2)$. Hence $K$ must lie in the closed
hemisphere of $ES$ centered at $x$.

If the point $x$ were to lie inside $N(K, \pi/2)$, then, just as in the previous section, we would find 
a contractible open subset $U$ of $S^3$ which separates $K$ from $K'$, which once again by 
Lemma E would imply that $\Hopf(f) = 0$.

Thus $x$ cannot lie inside $N(K, \pi/2)$, and it follows that no point of $K$ can lie in the open 
hemisphere of $ES$ centered at $x$. Since $K$ lies in that closed hemisphere, it must in fact lie 
on its boundary great circle $EC$.

We now assert that $K$ can not be a proper subset of $EC$, and see this in three cases 
as follows, supported by Figure 4 below.

\vfill
\eject

Assume for the moment that $K$ is a proper subset of the great circle $EC$.

\noindent
\bf Case 1\rm.\ \  $K'$ is disjoint from $EC$. Then an open 3-cell $U$ as shown in Figure 4 
separates $K$ from $K'$.
\begin{figure}[h!]
\center{\includegraphics[height=150pt]{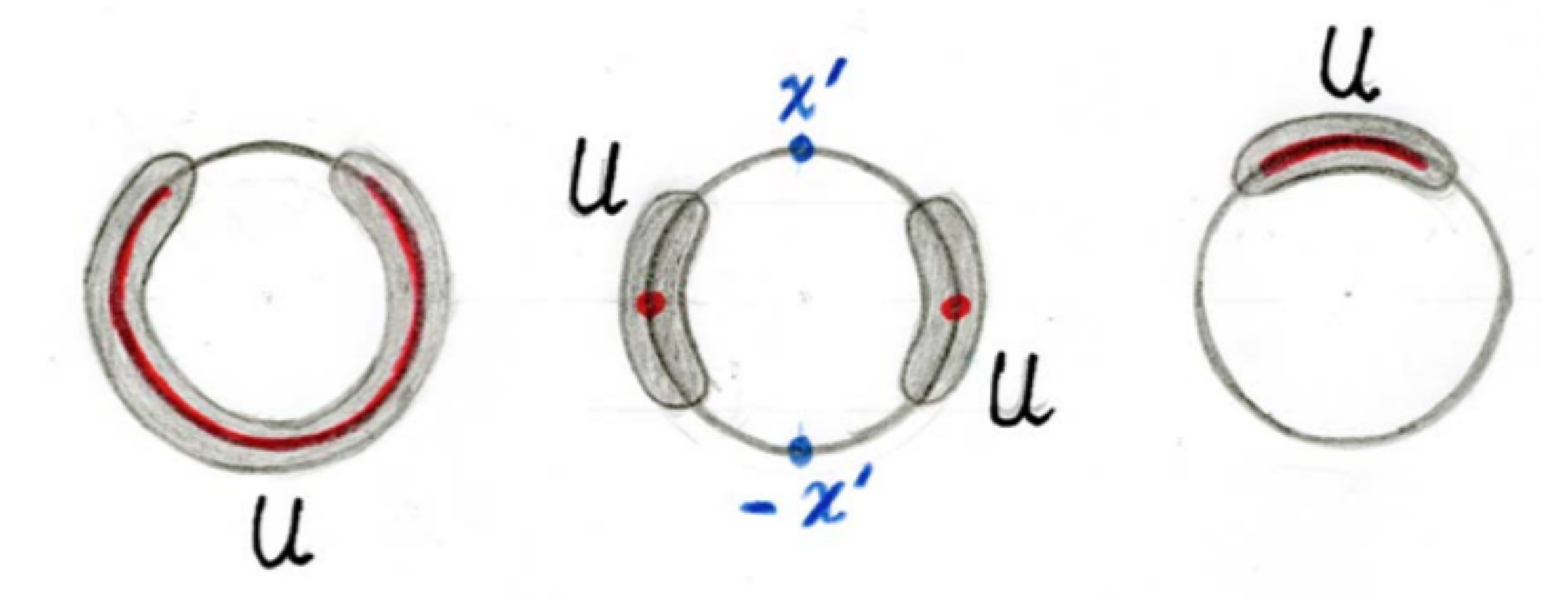}
\caption{\large{\textit{\textbf{\boldmath The fibre $K$ can not be a proper subset of the great circle $EC$}}}
}}
\end{figure}

\noindent
\bf Case 2\rm.\ \  $K'$ meets $EC$ in two antipodal points $x'$ and $-x'$. Then each of $x'$ 
and $-x'$ sits at the center of an open semi-circle on $EC$ which is forbidden to contain any points 
of $K$, since the distance in $S$ between $K$ and $K'$ is $\ge\pi/2$. So $K$ consists at most of 
two points, and then the disjoint union $U$ of two open three-cells, as shown in the figure, 
separates $K$ from $K'$.

\noindent
\bf Case 3\rm.\ \  $K'$ meets $EC$, but not just in two antipodal points. Then, as in Case 2 above, 
$K$ is forbidden to lie in a union of open semi-circles on $EC$, which in the present case is an open 
arc on $EC$. Therefore $K$ is constrained to lie in the complementary closed arc, and then the open 
3-cell $U$ shown in the figure separates $K$ from $K'$.

In each of the three cases above we have $H_1(U; \ints) = 0$, and then Lemma E would imply that 
$\Hopf(f) = 0$.

This contradiction shows that $K = f^{-1}(y)$ must be the entire great circle $EC$.

Since $y$ was an arbitrary point of $S^2$, we now know that all the fibres of $f$ are great 
circles in $S^3$.

\medskip

\noindent\bf Step 3. Any two great circle fibres of $f$ are parallel to one another.\rm

We claim now that any two great circle fibres of $f$ are \textit{{parallel}}, meaning that they are a constant 
distance apart from one another, and see this as follows.

Refer again to any pair of antipodal points $y$ and $-y$ on $S^2(1/2)$, and to their inverse images 
$K = f^{-1}(y)$ and $K' = f^{-1}(-y)$ in $S^3$, now known to be great circles.

No point of $K$ can be closer than $\pi/2$ to any point of $K'$, since their images $y$ and $-y$ 
under $f$ are exactly $\pi/2$ apart on $S^2(1/2)$ and we have $\Lip f\le 1$.

Thus the great circles $K$ and $K'$ on $S^3$ are orthogonal, meaning that they are the unit 
circles on a pair of orthogonal 2-planes through the origin in $\reals^4$.

Now let $z$ be a point on $S^2(1/2)$ at distance $\alpha$ from $y$ and at distance $\pi/2-\alpha$  
from $-y$.

Let $K'' = f^{-1}(z)$ be the corresponding great circle fibre. Where does $K''$ lie in $S^3$ with 
reference to $K$ and $K'$?

To answer that, consider the tubular neighborhoods $N(K,\alpha)$ and $N(K', \pi/2-\alpha)$ 
about $K$ and $K'$ in $S^3$. Each is an open solid torus, and their common boundary, call it 
$T_\alpha$, is a 2-dimensional torus, as shown in the figure below.

\begin{figure}[h!]
\center{\includegraphics[height=220pt]{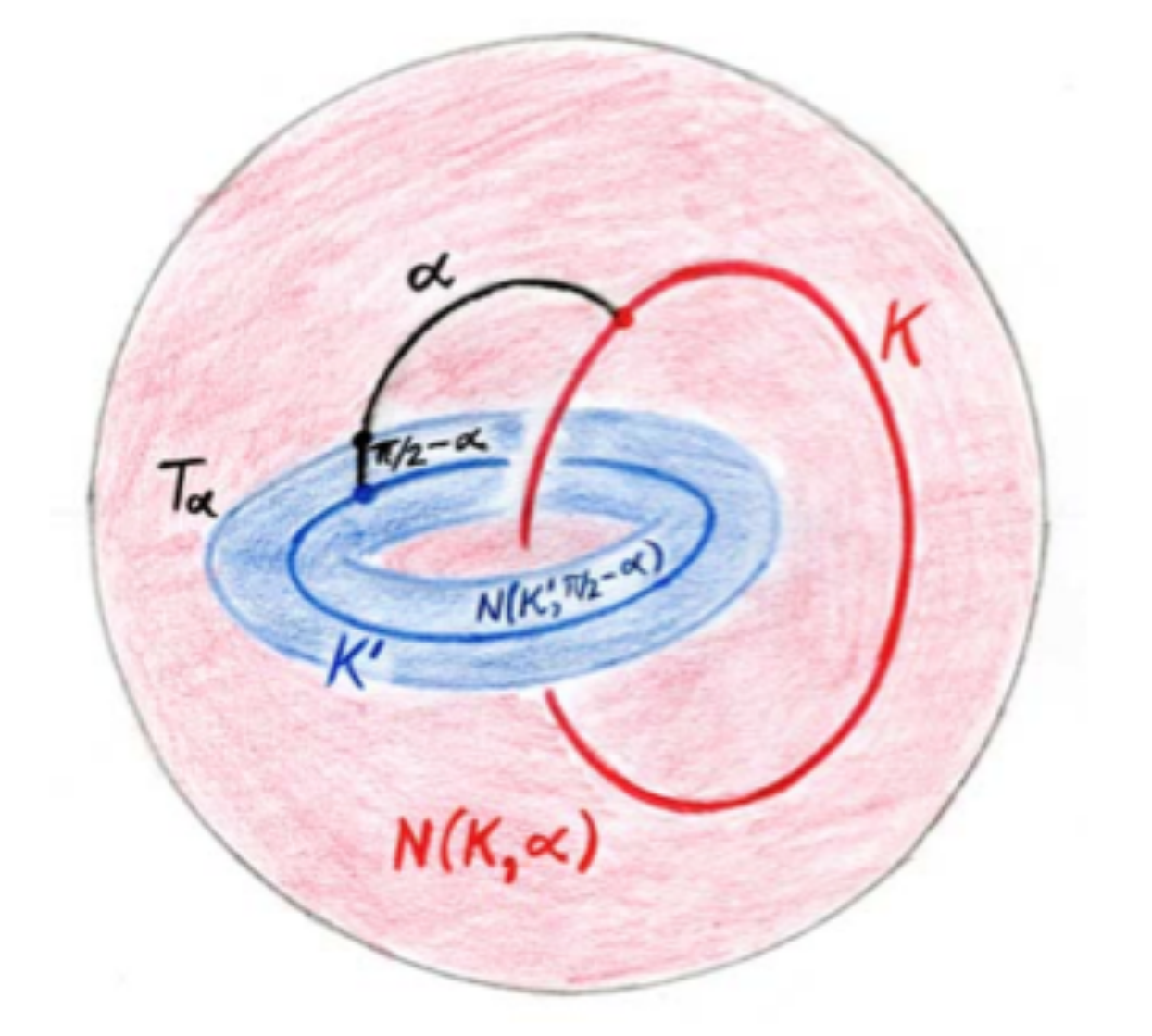}
\caption{\large{\textit{\textbf{\boldmath $S^3=N(K,\alpha)\cup T_\alpha\cup N(K',\pi/2-\alpha)$}}}}
}
\end{figure}

We claim that the great circle fibre $K'' = f^{-1}(z)$ must lie entirely on the \linebreak
2-dimensional torus 
$T_\alpha$.

If a part of $K''$ intrudes into the open set $N(K,\alpha)$, then that part is closer than $\alpha$ to 
$K$, yet is mapped by $f$ to the point $z$ which is exactly at distance $\alpha$ from $f(K) = y$, 
in contradiction to the assumption that $\Lip f \le 1$.

We get a similar contradiction if $K''$ intrudes into the open set $N(K', \pi/2-\alpha)$.

\vfill
\eject

Hence $K''$ lies on the common boundary $T_\alpha$ of these two open sets, 
and so is at constant distance $\alpha$ from $K$ and at constant distance $\pi/2-\alpha$ from $K'$.

Since $y$ and $z$ were arbitrary non-antipodal points on $S^2(1/2)$, we conclude that all the 
great circle fibres of the map $f$ are parallel to one another, as claimed.

\medskip

\noindent\bf Step 4. The map f is isometric to the Hopf projection.\rm

We have been considering a homotopically non-trivial map $f\colon S^3\to S^2(1/2)$ 
with Lipschitz constant $\le 1$, and have so far shown that the fibres of $f$ are parallel 
great circles on $S^3$.

But, as mentioned earlier, any fibration of $S^3$ by parallel great circles is isometric to the 
Hopf fibration. Thus $f$ induces a homotopically nontrivial map \linebreak
$\bah{f}\colon S^2(1/2)\to S^2(1/2)$ 
with Lipschitz constant $\le 1$, and this is easily seen to be an isometry. It follows that the map $f$ 
must be isometric to the Hopf projection.

This completes the proof of Theorem A in this first instance, and displays the style 
of argument that we will emulate for the general case.

\vfill
\eject

\centerline{\bf PART II. PROOF OF THEOREM A FOR ALL HOPF PROJECTIONS.}

\medskip

\noindent\bf The Hopf invariant\rm.

Look once again at our display of all the Hopf fibrations of round spheres by parallel great subspheres:
$$\begin{array}{lll}
\vphantom{\displaystyle{\int}}S^1 \subset S^3\to S^2=\cplxs\P^1,& S^1\subset S^5\ \to\cplxs\P^2,\ \ldots\,,&S^1\subset
S^{2n+1}\to\cplxs \P^n,\ldots\\
S^3 \subset S^7\to S^4=\quats \P^1,& S^3\subset S^{11}\to\quats \P^2,\ \ldots\,,&S^3
\subset S^{4n+3}\to \quats \P^n,\ldots\\
\vphantom{\displaystyle{\int}}S^7 \subset S^{15}\to S^8.& & 
\end{array}$$
We suppose that $f$ is a continuous map from 
$$S^{2n+1}\to\cplxs \P^n\quad \mbox{\rm or}\quad S^{4n+1}\to\quats \P^n\quad\mbox{\rm or}
\quad S^{15}	\to S^8,$$
and intend to give two equivalent definitions of its Hopf invariant. 

Before doing that, we recall the topology of the base spaces.

A choice of ascending complex vector spaces $\cplxs^1\subset\cplxs^2\subset\cplxs^3\subset\cdots$
leads to an ascending sequence of complex projective spaces 
$\cplxs \P^0 \subset \cplxs \P^1\subset \cplxs \P^2 \subset \cdots$. 
The cohomology ring $H^\dast(\cplxs \P^n; \ints)$ is a 
truncated polynomial ring with one generator $[\omega]$ in dimension 2 and with the relation 
$[\omega]^{n+1}=0$. We can take $\omega$ to be the K\"ahler form, scaled so that it integrates to 
1 over $\cplxs \P^1 = S^2(1/2)$.

A choice of ascending quaternionic vector spaces 
$\quats^1\subset\quats^2\subset\quats^3\subset\cdots$ leads to an ascending sequence of 
quaternionic projective spaces $\quats \P^0\subset \quats \P^1\subset \quats \P^2\subset\cdots$. 
The cohomology ring $H^\dast(\quats \P^n; \ints)$ is a truncated polynomial ring with one 
generator $[\omega]$ in dimension 4 and with the relation $[\omega]^{n+1}= 0$. 
We can take $\omega$ to be the quaternionic K\"ahler form, scaled so that it integrates to 1 
over $\quats \P^1	= S^4(1/2)$.

For simplicity of expression and to gain the advantage of making our arguments
more concrete, we will focus on maps 
$f\colon S^{2n+1}	\to\cplxs \P^n$, and then comment afterwards on the very slight changes 
needed to handle maps of $S^{4n+1}\to\quats \P^n$ and of $S^{15}\to S^8$.

\noindent \bf (1) \rm\ \  Given a map $f \colon S^{2n+1}\to \cplxs \P^n$, homotope it so that it is 
smooth, has a given value $y$ in $\cplxs \P^n$ as regular value, and so that it is transverse to 
the corresponding ``antipodal'' $\cplxs \P^{n-1}$, which is simply the cut locus of $y$ in the 
usual Riemannian metric on $\cplxs \P^n$. Then the inverse image $K = f^{-1}(y)$ is a finite union 
of smooth simple closed curves in $S^{2n+1}$, while the inverse image $K' = f^{-1}(\cplxs \P^{n-1})$ 
is a smooth submanifold of $S^{2n+1}$ of dimension $2n-1$. 

\vfill
\eject

Orienting $S^{2n+1}$ arbitrarily, and
$\cplxs \P^n$ in the usual way, we derive orientations for $K$ and $K'$ just as we did for maps of 
$S^3\to S^2$. Then the Hopf invariant of $f$ is defined to be the total linking number of all 
the components of $K$ with all the components of $K'$.

This definition is independent of the choice of $y$ and $\cplxs \P^{n-1}$ in $\cplxs \P^n$, and 
depends only on the homotopy class of $f$.

\noindent\bf (2) \rm\ \  Let $\omega$ be the K\"ahler form on $\cplxs \P^n$, scaled so that 
$\ip{\omega}{\cplxs \P^1}=1$. Then the pullback $f^\dast\omega$ is a 2-cocycle on $S^{2n+1}$. 
Since $H^2(S^{2n+1};\ints) = 0$, there is a 1-dimensional integral cochain 
$\alpha$ on $S^{2n+1}$ such that $d\alpha=f^\dast\omega$.
Then the integer $\ip{\alpha\cup (f^\dast\omega)^n}{S^{2n+1}}$ is defined to be the Hopf invariant of $f$.

One shows that this definition is independent of the choice of 2-cocyle $\omega$ on $\cplxs \P^n$	
which generates $H^2(\cplxs \P^n; \ints) \cong \ints$, and of the choice of 1-cochain $\alpha$ on 
$S^{2n+1}$ such that $d\alpha = f^\dast\omega$, and that it depends only on the homotopy class
of $f$.

If f is smooth, we can use de Rham cohomology for this approach, just as we did for maps of 
$S^3\to S^2$.

\medskip

\noindent\bf Mix-and-match formula for the Hopf invariant.\rm

The situation here is the same as for maps of $S^3\to S^2$. 

Suppose $\omega'$ is another 2-cocycle on $\cplxs \P^n$ with $\ip{\omega'}{\cplxs \P^1}=1$.

Then instead of the above formula
$$\Hopf (f) = \ip{\alpha\cup(f^\dast\omega)^n}{S^{2n+1}}$$
for the Hopf invariant, we have the mix-and-match formula,
$$\Hopf(f) = \ip{\alpha\cup(f^\dast\omega')^n}{S^{2n+1}}.$$
To verify this, first write $\omega-\omega'=d\eta$ for some 1-cochain $\eta$ on $\cplxs \P^n$.

It follows that $\omega^n-(\omega')^n=d\zeta$ for some $(2n-1)$-cochain $\zeta$ on 
$\cplxs \P^n$. 

Write
$$\alpha\cup(f^\dast\omega)^n-\alpha\cup(f^\dast\omega')^n=\alpha\cup f^\dast d\zeta
=\alpha\cup d\,f^\dast\zeta,$$
and then the integration by parts given earlier in the case of $S^3\to S^2$, now with $\zeta$ in place 
of $\eta$ there, finishes the present argument and confirms the mix-and-match formula above.

\vfill
\eject

\noindent\bf A sufficient condition for the Hopf invariant to be zero. 

\medskip

\noindent Lemma E (Complete version).

\noindent (1) \it\ \ Let $f\colon S^{2n+1}\to\cplxs \P^n$ be a continuous map, and let $y$ and 
$\cplxs \P^{n-1}$ be a point and disjoint projective hyperplane in $\cplxs \P^n$, 
with inverse images $K = f^{-1}(y)$ and $K' = f^{-1}(\cplxs \P^{n-1})$. 
Suppose there is an open set $U$ in $S^{2n+1}$ such that
$$K\subset U \subset \bah{U}\subset S^{2n+1}-K' \quad \mbox{\rm and}\quad  
H_1(U;\ints)=0.$$ 
Then the Hopf invariant of $f$ is zero.

\noindent\bf (2) \it\ \  Let $f\colon S^{4n+1}\to\quats \P^n$ be a continuous map, and let $y$ and 
$\quats \P^{n-1}$ be a point and disjoint projective hyperplane in $\quats \P^n$, 
with inverse images $K = f^{-1}(y)$ and $K' = f^{-1}(\quats \P^{n-1})$. 
Suppose there is an open set $U$ in $S^{4n+1}$ such that
$$K\subset U\subset\bah{U}\subset S^{4n+1}- K' \quad\mbox{\rm and}\quad 
H_3(U;\ints)=0.$$ 
Then the Hopf invariant of $f$ is zero.

\noindent\bf(3)\it\ \  Let $f\colon S^{15}\to S^8$ be a continuous map, and let $y$ and $y'$ be two 
points of $S^8$, with inverse images $K = f^{-1}(y)$ and $K' = f^{-1}(y')$. Suppose there is an 
open set $U$ in $S^{15}$ such that
$$K\subset U\subset\bah{U} \subset S^{15}-K' \quad\mbox{\rm and}\quad H_7(U;\ints)=0.$$ 
Then the Hopf invariant of $f$ is zero. \rm

The proof is the same as for the prototype discussed earlier.

\medskip

\noindent\bf Beginning the proof of Theorem A.\rm

We will give the proof for maps of $S^5\to\cplxs \P^2$, leaning heavily on the techniques developed 
for the case $S^3\to S^2$, and afterwards explain the small adjustments needed to handle 
the general case.

We begin with a continuous map $f\colon S^5\to\cplxs \P^2$ with nonzero Hopf invariant, 
assume that $\Lip f\le 1$ and set out to prove that $\Lip f = 1$ and that $f$ is isometric to the 
Hopf projection.

\medskip

\noindent\bf
Step 1. Each fibre of $f$ lies on a great 4-sphere in $S^5$.\rm

In $\cplxs \P^2$, we focus on an arbitrary point $y$ and on its cut locus $Y'=\cplxs \P^1\cong S^2$ 
at maximal constant distance $\pi/2$ along every geodesic streaming out from $y$.

In $S^5$ we focus on the fibre $K = f^{-1}(y)$ and on the union of fibres $K' = f^{-1}(Y')$.

Let $N(K, \pi/2)$ again denote the open $\pi/2$ neighborhood of $K$ in $S^5$.

Since $\Lip f \le 1$, no point in $N(K,\pi/2)$ can map to $Y'$.

\medskip

On the other hand, some point in $S^5$ must map to $Y'$, because otherwise the image of $f$ 
would lie in $\cplxs \P^2 -Y'$, which is an open 4-cell, and this would make $f$ homotopically trivial.

Say $f(-x)\in Y'$.

Since $\Lip f\le 1$, the point $-x$ can not lie in $N(K, \pi/2)$, and therefore no point of $K$ can lie in 
the open hemisphere $N(-x, \pi/2)$. Hence $K$ must lie in the closed hemisphere of 
$S^5$ centered at $x$.

If the point $x$ were to lie in $N(K, \pi/2)$ then, just as in the case of maps from \linebreak
$S^3\to S^2$, we could construct a contractible open neighborhood $U$ of $K$ which separates it 
from $K'$, and then conclude from Lemma E that the Hopf invariant of $f$ must be zero.

Thus the point $x$ can not lie in $N(K, \pi/2)$, and it follows that no point of $K$ can lie in the open hemisphere $N(x, \pi/2)$. Since we already know that $K$ lies in that closed hemisphere, it follows 
that $K$ must lie on its boundary, which is a great 4-sphere $ES^4$ in $S^5$.

\medskip

\noindent\bf Steps 2, 3, 4. Each fibre of $f$ is a great circle in $S^5$.\rm

We then follow the argument from the case of maps from $S^3\to S^2$, using the fact that a small 
open neighborhood of $ES^4$ has trivial 1-dimensional homology, and invoke Lemma E once again 
to conclude that $K$ must in fact lie on a great 3-sphere $ES^3\subset ES^4$.

We iterate this twice more to conclude that K must lie on a great circle $ES^1$, and then copy our 
earlier argument from the $S^3\to S^2$ case to conclude that $K$ can not be a proper subset of 
that great circle, and hence must coincide with it.

Since $y$ was an arbitrary point of $\cplxs \P^2$, we now know that each fibre 
$K = f^{-1}(y)$ is a great circle on $S^5$.

\vfill
\eject

\noindent\bf Step 5. Any two great circle fibres of $f$ are parallel to one another.\rm

Consider the great circle $K = f^{-1}(y)$ and the set $K' = f^{-1}(Y')$, which must be a union of 
great circles.
Since $\Lip f \le 1$, the set $K'$ must lie within the great 3-sphere $S^3$ in $S^5$ which is orthogonal 
to $K$ and at constant maximal distance $\pi/2$ from it.

If $K'$ were a proper subset of $S^3$, we could easily construct a contractible open set $U$ in 
$S^5$ which separates $K'$ from $K$, and then conclude from Lemma E that the Hopf invariant 
of $f$ must be zero.

Hence $K' = S^3$.

We now copy the argument from the $S^3\to S^2$ case to conclude that any two great circle fibres 
of $f$ are parallel to one another.

\medskip

\noindent\bf Step 6. The map $f$ is isometric to the Hopf projection.\rm

Just as in the $S^3\to S^2$ case, this follows from the known fact, mentioned earlier, that any 
fibration of a round sphere by parallel great subspheres is isometric to the corresponding Hopf 
fibration.

\medskip

\noindent\bf Completion of the proof of Theorem A.\rm

The same argument handles all the Hopf fibrations, and in each case shows that a map
$$S^{2n+1}\to\cplxs \P^n\quad\mbox{\rm or}\quad S^{4n+3}\to\quats \P^n\quad
\mbox{\rm or}\quad S^{15}\to S^8$$
with nonzero Hopf invariant and Lipschitz constant $\le 1$ must have Lipschitz constant 
equal to 1 and be isometric to the corresponding Hopf projection.

This completes the proof of Theorem A.

\medskip

\noindent
\textit{\textbf{Comment.}}\rm

The set of homotopy classes of maps from $S^{2n+1}\to\cplxs \P^n$ is in one-to-one correspondence 
with the integers, as one sees from the homotopy sequence of the bundle 
$S^1\subset S^{2n+1}	\to\cplxs \P^n$, with the Hopf invariant providing the correspondence.

But in the remaining cases, there are homotopically nontrivial maps which nevertheless have 
zero Hopf invariant.

Consider for example the Hopf fibration $S^3\subset S^7\to S^4$. From the homotopy sequence of 
this bundle and the fact that the fibre is contractible in the total space, we get
$$[S^7,S^4]\,\cong \,\pi_7(S^4)\,\cong\,\pi_7(S^7) + \pi_6(S^3)\,\cong \,\ints + \ints_{12}.$$ 
The $\ints$-summand of $\pi_7(S^4)$ corresponds to the Hopf invariant, but 
the maps in the $\ints_{12}$-summand all have Hopf invariant zero. And likewise for all the 
quaternionic Hopf projections. In the one remaining case, $S^7\subset S^{15}\to S^8$, we get
$$[S^{15}, S^8]\,\cong \,\pi_{15}(S^8)\,\cong\,\pi_{15}(S^{15}) + \pi_{14}(S^7)\,\cong\,
\ints + \ints_{120},$$ 
with the $\ints$-summand corresponding to the Hopf invariant, but with all 
the maps in the \linebreak $\ints_{120}$-summand having Hopf invariant zero.

\vfill
\eject

\centerline{\bf PART III. PROOF OF THEOREM B.\rm}

\medskip

\noindent\bf Statement of the Key Lemma.\rm

Theorem B metrically characterizes the inclusion map $i\colon \Delta S^n\to S^n\times S^n$ 
of the diagonal as a Lipschitz constant minimizer in its homotopy class, unique up to composition 
with isometries of domain and range.

To prove this, we start with a map $f\colon \Delta S^n \to S^n\times S^n$ which is homotopic to $i$, assume that $\Lip f \le 1$, and aim to show that $\Lip f = 1$ and that $f$ is isometric to $i$.

The basic tool is the Key Lemma, stated below. 

If $x$ is any point on $S^n$, then $-x$ is the antipodal point, and their distance
apart on $S^n$ is $\pi$.

Likewise, if $(x, y)$ is any point of $S^n\times S^n$, then $(-x ,-y)$ will be called its 
\it antipodal point\rm, and their distance apart on $S^n\times S^n$ is $\pi\sqrt{2}$. 
This is the maximum distance between any two points of $S^n\times S^n$.

\medskip

\noindent\bf Key Lemma.\ \  \textit{\textbf{\boldmath Let $f\colon\Delta S^n\to S^n\times S^n$ 
be a map which is homotopic 
to the inclusion. Then its image $f(\Delta S^n)$ contains a pair of antipodal points $(x,y)$ and 
$(-x,-y)$ in $S^n\times S^n$.}} \rm

\medskip

The claim, in other words, is that the image $f(\Delta S^n)$ contains a pair of points at 
maximum distance apart in $S^n\times S^n$. In applying the Key Lemma, these two points will 
serve as a kind of framework, upon which the image is stretched.

The Key Lemma has a Borsuk-Ulam flavor.

\medskip

\noindent\bf Suspension.\rm

Consider a map $\ph\colon S^m\to S^n$. Then a concrete model for the suspension 
of $\ph$ is the map $\Sigma\ph\colon S^{m+1}\to S^{n+1}$ defined by
$$\Sigma\ph\,(x\cos t,\sin t)=(\ph(x)\cos t,\sin t),$$
where $x\in S^m$ and $-\pi/2\le t\le \pi/2$, as illustrated in the figure below.

\begin{figure}[h!]
\center{\includegraphics[height=145pt]{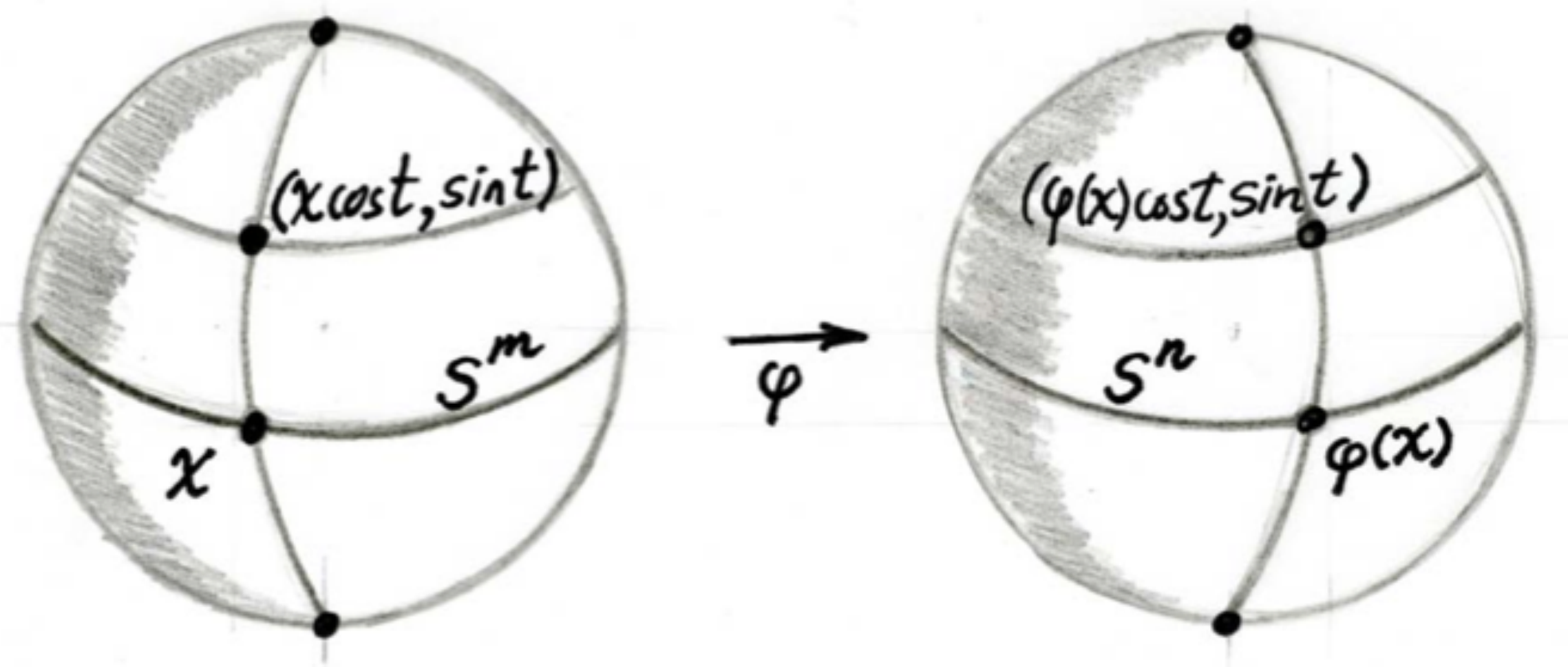}
\caption{\large{\textit{\textbf{\boldmath $\Sigma\ph\colon S^{m+1}\to S^{n+1}$}}}}
}
\end{figure}

The suspension of  $\ph$ takes $m$-spheres of constant latitude on $S^{m+1}$ to $n$-spheres of
constant latitude on $S^{n+1}$ by rescaled copies of $\ph$. It is almost never smooth at the 
north and south poles, no matter how smooth $\ph$ is.

\vfill
\eject

If $f = (f_1 ,f_2)\colon S^k\to S^m\times S^n$, then we define the suspension of $f$ to be the map
$$\Sigma f = (\Sigma f_1,\Sigma f_2)\colon S^{k+1}\to S^{m+1}\times  S^{n+1}.$$
If $f$ = $(f_1 , f_2)\colon  \Delta S^n\to S^n\times S^n$ is homotopic to the inclusion, then \linebreak
$\Sigma f\colon \Delta S^{n+1} \to S^{n+1}\times S^{n+1}$ is also homotopic to the inclusion, 
with the obvious rescaling to make the suspension of $\Delta S^n$ into $\Delta S^{n+1}$.

\medskip

\noindent\bf Plan of the proof of Theorem B.\rm

We start with a map $f\colon \Delta S^n\to S^n\times S^n$ which is homotopic to the inclusion,
assume that $\Lip f\le 1$, and aim to show that $\Lip f = 1$ and that $f$ is isometric to the inclusion.

The argument is by induction on $n$.

We assume the truth of the Key Lemma, leave the base step $n = 1$ as an exercise for the reader, 
and begin with the induction step as follows.

For $n > 1$, we use the Key Lemma together with the hypothesis that $\Lip f\le 1$ to
\it desuspend \rm $f$ to a map
$f'\colon \Delta S^{n-1}\to S^{n-1}\times S^{n-1}$, meaning that 
$\Sigma f' = f$, such that $f'$ is homotopic to the inclusion and
satisfies $\Lip f' \le 1$. 

Then by the induction hypothesis, we know that $f'$ is isometric to 
the inclusion, and immediately conclude the same for $f =\Sigma f'$. 

Finally, we give the proof of the Key Lemma.

For even $n$, this is a straightforward intersection argument in $S^n\times S^n$ using\linebreak homology 
with integer coefficients.

For odd $n$, this is an intersection argument in the symmetric product $S^n \dst S^n$ using 
homology with coefficients mod 2 .

\medskip

\noindent\bf The induction step.\rm

We assume the truth of Theorem B for $n-1$, and show how to prove it for $n$. 

We start with a map $f\colon \Delta S^n\to S^n\times S^n$ which is homotopic to the inclusion
and satisfies $\Lip f \le 1$. 

By the Key Lemma, the image $f(\Delta S^n)$ contains a pair of 
antipodal points $(x,y)$ and $(-x,-y)$ in $S^n\times S^n$. 

Their distance apart in 
$S^n\times S^n$	is $\pi\sqrt{2}$, and since $\Lip f \le 1$, they must be
the images of a pair of antipodal points, say $u$ and $-u$, in $\Delta S^n$.

\begin{figure}[h!]
\center{\includegraphics[height=170pt]{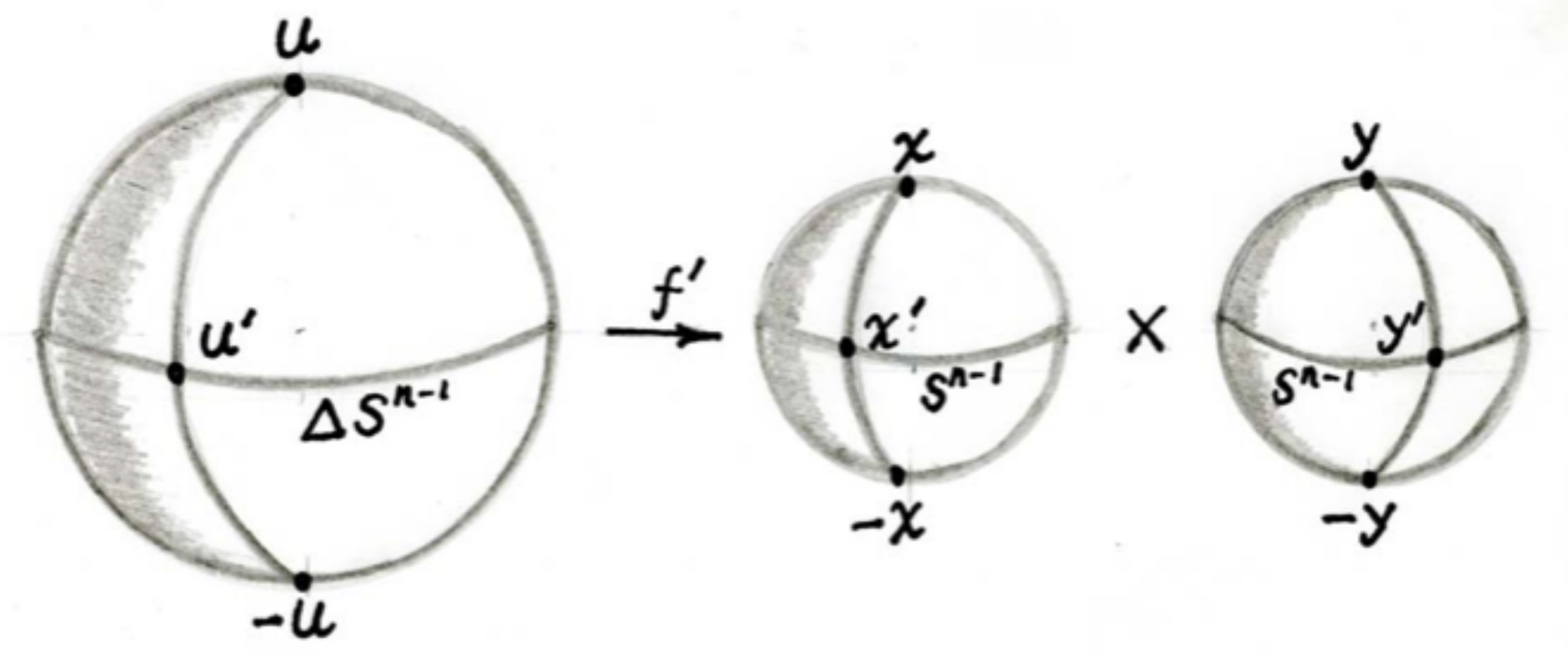}
\caption{\large{\textit{\textbf{\boldmath $f\colon\Delta S^n\to S^n\times S^n$}}}
}}
\end{figure}

On $\Delta S^n$, each semicircle from $u$ to $-u$ is a geodesic of length $\pi\sqrt{2}$, 
and since $\Lip f \le 1$, it must be taken by $f$ to a geodesic, also of length $\pi\sqrt{2}$, 
from $(x,y)$ to $(-x,-y)$ on $S^n\times S^n$.

The first coordinate of this image geodesic on $S^n\times S^n$ runs from $x$ to $-x$ along a 
semicircle on the first $S^n$ factor, and likewise the second coordinate runs from $y$ 
to $-y$ along a semicircle on the second $S^n$ factor.

Since $\Lip f \le 1$, the map $f$ from the semicircle on $\Delta S^n$ to the product of the two 
semicircles on $S^n\times S^n$ must be distance-preserving, with no leeway for slowing down or speeding up.

Let $\Delta S^{n-1}$ denote the equatorial $(n-1)$-sphere on $\Delta S^n$ with poles at $u$ and 
$-u$, and likewise let $S^{n-1}$ denote (ambiguously) the equatorial $(n-1)$-spheres on the two 
$S^n$ factors, with poles at $x$ and $-x$, and at $y$ and $-y$, respectively.

Let $u'$, $x'$ and $y'$ denote the points where the three semicircles meet their respective 
equators, as shown in the figure above.
Then $f(u') = (x', y')$. So we define 
$$f'\colon \Delta S^{n-1}\to S^{n-1}\times  S^{n-1}$$
to be the restriction of $f$ to the equator $\Delta S^{n-1}$ on $\Delta S^n$. 

Then we see from the above construction that $f$ is the suspension of $f'$,
that is, $f = \Sigma f'$. 

Since $\Delta S^{n-1}$ is totally geodesic in $\Delta S^n$, and $S^{n-1}\times S^{n-1}$ 
is totally geodesic  in $S^n\times S^n$, the hypothesis that $\Lip f\le 1$ implies that 
$\Lip f'\le 1$.

Furthermore, the hypothesis that $f\colon \Delta S^n\to S^n\times S^n$ is homotopic to the 
inclusion implies that $f'\colon \Delta S^{n-1}\to S^{n-1}\times S^{n-1}$ is also homotopic 
to the inclusion.

The induction hypothesis, that Theorem B is true in dimension $n-1$, now tells us that $f'$ must be isometric to the inclusion 
$i'\colon \Delta S^{n-1}\to S^{n-1}\times S^{n-1}$, and it follows immediately that $f =\Sigma f'$ must 
be isometric to the inclusion $i\colon \Delta S^n\to  S^n\times S^n$.

This completes the proof of Theorem B, modulo the Key Lemma.

\medskip

\noindent\bf Proof of the Key Lemma for even $n$.\rm

We start with a map $f\colon \Delta S^n\to S^n\times S^n$ which is homotopic to the inclusion, 
and must find a pair of antipodal points $(x, y)$ and $(-x, -y)$ in its image.

Let $a$ and $b$ denote the generators of $H_n(S^n\times S^n; \ints)$ represented by 
$S^n \times \mbox{\it point}$ and by $\mbox{\it point}\times S^n$. Then $f(\Delta S^n)$ 
can be regarded as a singular $n$-cycle representing the class $a$ + $b$.

\begin{figure}[h!]
\center{\includegraphics[height=140pt]{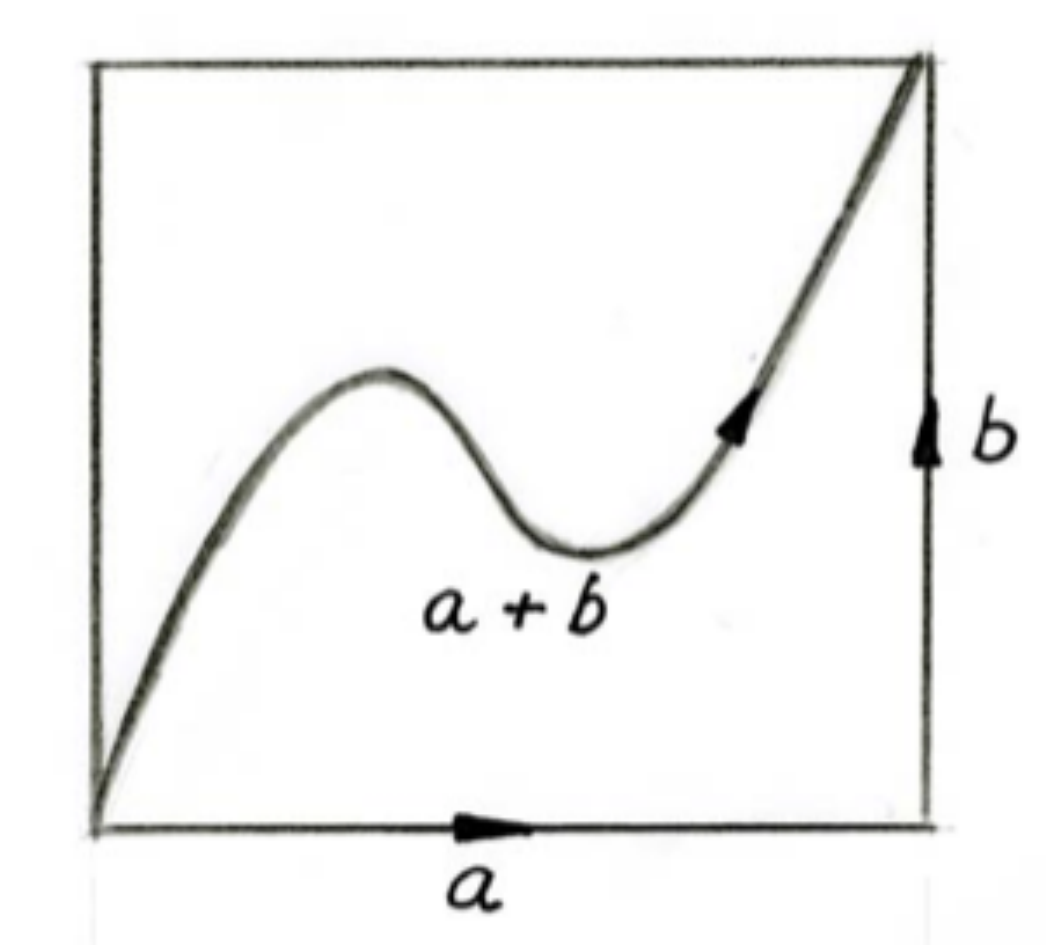}
\caption{\large{\textit{\textbf{The $n$-cycle $f(\Delta S^n)$ represents the homology class $a+b$}}}
}}
\end{figure}

Since $n$ is even, the intersection form on $H_n(S^n\times S^n; \ints)$ is given by
$$a\cdt a=b\cdt b=0,\quad a\cdt b=b\cdt a=1.$$ 
If $f = (f_1 ,f_2)$, let $-f = (-f_1 ,-f_2)$, so that
$$-f(\Delta S^n) = \{(-x, -y)\,:\, (x, y) \in f (\Delta S^n)\}.$$ 
In other words, $-f(\Delta S^n)$ consists of the antipodes of all the points in $f(\Delta S^n)$.

Since $n$ is even, the oriented singular $n$-cycle $-f(\Delta S^n)$ represents the 
homology class $-a-b$.

The intersection number of the singular $n$-cycles $f(\Delta S^n)$ and $-f(\Delta S^n)$ is 
$$(a+b)\cdt(-a-b) = -2.$$
Hence $f(\Delta S^n)$ and $-f(\Delta S^n)$ certainly have a nonempty intersection. 

Therefore, for some point $(x, y)$ in $f(\Delta S^n)$, the point $(-x, -y)$ is also in
$f(\Delta S^n)$, which is precisely the claim of the Key Lemma. 

\medskip

\noindent\bf Comment\rm.\ \  When $n$ is odd, the intersection form on 
$H_n(S^n\times S^n; \ints)$ is given by
$$a\cdt a=b\cdt b=0,\quad a\cdt b=1,\quad b\cdt a=-1,$$ 
and the singular $n$-cycles $f(\Delta S^n)$ and $-f(\Delta S^n)$ both represent the same
class $a+b$. 

The intersection number of these two $n$-cycles is therefore
$$(a+b)\cdt (a+b) = 0,$$ 
and the preceding argument falls apart.

\medskip

\noindent\bf Rephrasing the Key Lemma.\rm

First we restate it.

\medskip

\noindent\bf Key Lemma.\ \  \textit{\textbf{\boldmath Let $f\colon \Delta S^n\to  S^n \times S^n$ 
be a map which is 
homotopic to the inclusion. Then its image $f(\Delta S^n)$ contains a pair of antipodal points 
$(x,y)$ and $(-x,-y)$ in $S^n\times S^n$.}}\rm

\medskip

Then we rephrase it.

\vfill
\eject

\noindent\bf Key Lemma, rephrased.\ \  \textit{\textbf{\boldmath Let $f_1$ and 
$f_2\colon S^n\to S^n$ be two maps which 
are both homotopic to the identity. Then there are points $u$ and $v$ in $S^n$ such that 
$f_1(u)$ and $f_1(v)$ are antipodal, and at the same time $f_2(u)$ and $f_2(v)$ are also antipodal.}}
\rm

\medskip

To match the two versions, put $f(x, x) = (f_1(x), f_2(x))$.

In the two hypotheses, $f$ is homotopic to the inclusion if and only if $f_1$ and $f_2$ are both 
homotopic to the identity.

In the two conclusions, $f(\Delta S^n)$ contains a pair of antipodal points, call them 
$f(u, u) = (f_1(u), f_2(u))$ and $f(v, v) = (f_1(v), f_2(v))$, if and only if $f_1(u)$ and $f_1(v)$ 
are antipodal, and at the same time $f_2(u)$ and $f_2(v)$ are antipodal.

\medskip

The proof of the Key Lemma which we give here is due to Dennis Sullivan. It begins with the alternative
phrasing above, and then moves the scene of action from the cartesian product $S^n\times S^n$ down 
to the symmetric product $S^n\dst S^n$, in which every point $(u, v)$ is identified with its ``flip'' $(v, u)$.

The virtue of this move is that the image of the ``anti-diagonal'' in $S^n\dst S^n$ will be seen to have 
self-intersection number 1 mod 2 there, independent of the parity of $n$.

To prepare for the argument, we pause to discuss the geometries of both the cartesian and symmetric 
products.

\medskip

\noindent\bf Geometry of the cartesian product $S^n\times S^n$.\rm 

Let
$$\eqalign{
D &= \Delta S^n = \mbox{\it diagonal $n$-sphere} = \{(x,x)\,:\,x\in S^n\} \subset S^n\times S^n,\cr 
A &= \mbox{\it anti-diagonal $n$-sphere} = \{(x, -x)\,:\, x \in S^n\}\subset S^n\times S^n.
}$$
Each of $D$ and $A$ is the focal locus of the other in $S^n\times S^n$, and the 
isometry \linebreak
$\ph\colon S^n\times S^n\to S^n\times S^n$ defined by $\ph(x, y) = (x,-y)$ interchanges them.

The diagonal $D$ and anti-diagonal $A$ are homologous to one another when $n$ is odd, but not 
when $n$ is even.

The set $U = \{(x,y)\in S^n\times S^n\,:\, x\cdt y = 0\}$ is a copy of the Stiefel manifold 
$V_2\reals^{n+1}$ of orthonormal two-frames in $R^{n+1}$, and is situated halfway between 
$D$ and $A$ in $S^n\times S^n$.

\begin{figure}[h!]
\center{\includegraphics[height=140pt]{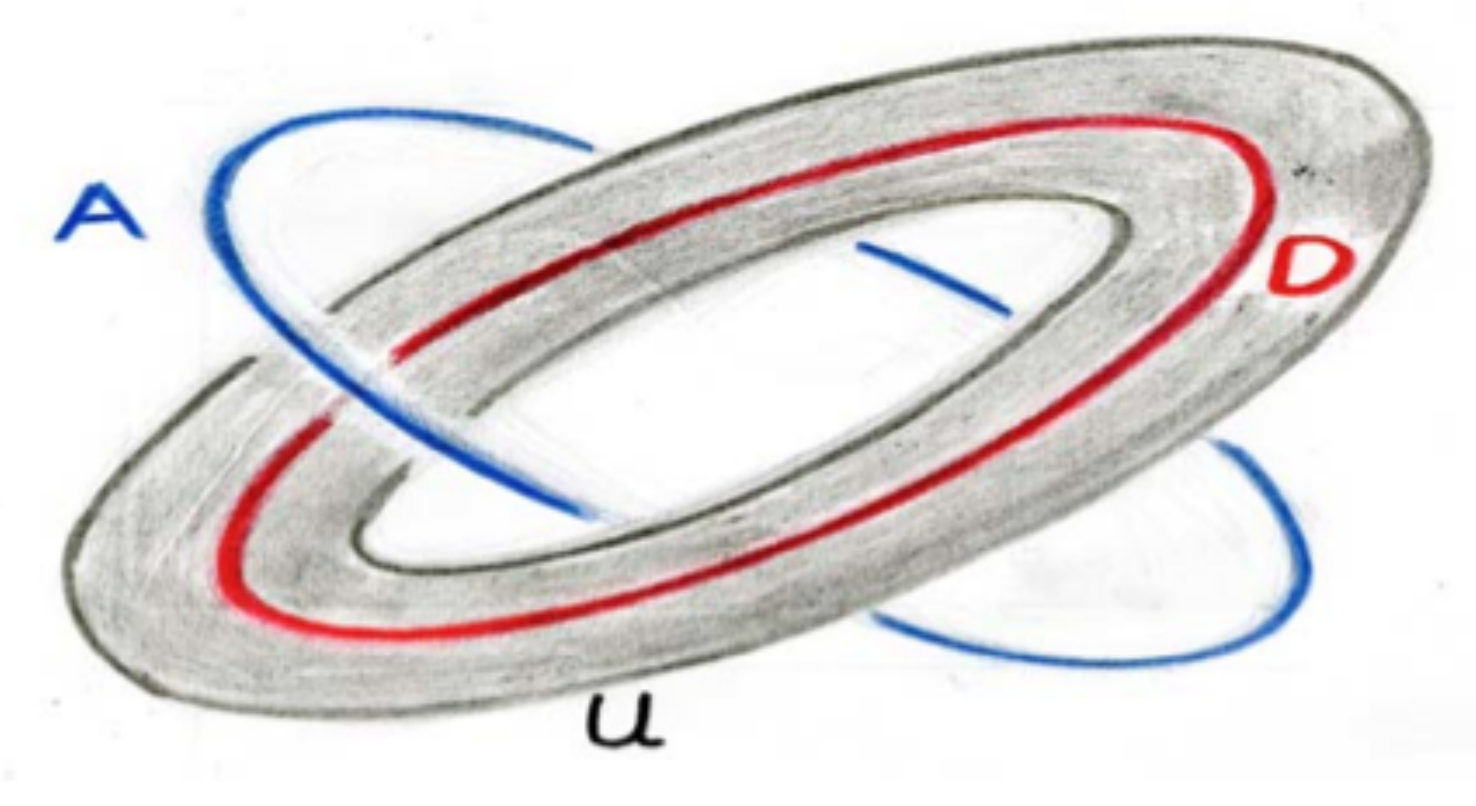}
\caption{\large{\textit{\textbf{\boldmath The hypersurface $U$ is halfway between $D$ and $A$ in 
$S^n\times S^n$}}}}
}
\end{figure}

The inner product function $\IP\colon S^n\times S^n\to\reals$ defined by $\IP(x, y) = x \cdt y$
takes values in the interval $[-1, 1]$, and we have 
$$D = \IP^{-1}(1)\quad U = \IP^{-1}(0), \quad A = \IP^{-1}(-1).$$
The level sets $\IP^{-1}(t)$ for $-1 < t < 1$ are all homeomorphic to one another, and together 
foliate the complement of $D$ and $A$ in $S^n\times S^n$, collapsing to $D$ at one end and to $A$
at the other.

The sets $N(D) = \IP^{-1}([0, 1])$ and $N(A) = \IP^{-1}([-1, 0])$ are closed
tubular neighborhoods of 
$D$ and $A$ in $S^n\times S^n$ which share $U$ as a common boundary. They are both copies of 
the unit disk bundle of the tangent bundle of $S^n$.

\medskip

\noindent\bf The symmetric product $S^n\dst S^n$.\rm

 Consider the involution $\sigma\colon S^n\times S^n\to S^n\times S^n$ defined by 
 $\sigma(x, y) = (y, x)$.
 
The \it symmetric product \rm $S^n \dst S^n$ is obtained from $S^n\times S^n$ by dividing out by 
this involution,
$$S^n \dst S^n = \frac{S^n \times S^n}{(x,y)\sim(y,x)}.$$ 
Let $p\colon S^n\times S^n \to S^n \dst S^n$ be the projection map, and write $p(x,y) = [x,y]$.

Each level set of the inner product function $\IP\colon S^n \times S^n\to [-1, 1]$ is invariant under 
the involution $\sigma$, and hence the decomposition of $S^n\times S^n$ into these level sets 
projects under $p$ to a corresponding decomposition of $S^n \dst S^n$.

The fixed point set of $\sigma$ is the diagonal $n$-sphere $D$, which projects one-to-one to its image 
$p(D) = D'$ in $S^n \dst S^n$. At the other extreme is the anti-diagonal $A$ in $S^n \times S^n$, 
which projects two-to-one to its image $p(A) = A'$, a copy of $\reals \P^n$, in $S^n \dst S^n$.

The tubular neighborhoods $N(D)$ and $N(A)$ in $S^n \times S^n$ project to tubular 
neighborhoods $N(D')$ and $N(A')$ in $S^n \dst S^n$.

The symmetric product $S^1 \dst S^1$ is a M\"obius band, with $D'$ as its boundary, while the 
symmetric product $S^2 \dst S^2$ is homeomorphic to $\cplxs \P^2$.

By contrast, the symmetric product $S^n \dst S^n$ fails to be a manifold along $D'$ starting with 
$n = 3$. Nevertheless, for all $n$, $(S^n\dst S^n) - D'$ is a (noncompact) manifold containing $A'$ 
as a submanifold.

Since the involution $\sigma$ of $S^n\times S^n$ is orientation-preserving for even $n$ and 
orientation-reversing for odd $n$, the symmetric product $S^n\dst  S^n$ is orientable for 
even $n$ and non-orientable for odd $n$.

By contrast, the image $A'$ of the anti-diagonal is homeomorphic to $\reals \P^n$, and is therefore 
non-orientable for even $n$ and orientable for odd $n$ .

Thus the symmetric product $S^n \dst  S^n$ has some prominent non-orientable feature for all $n$.

\medskip

\noindent\bf Self-intersection number of the anti-diagonal in $S^n \dst  S^n$.\rm

Recall that in $S^n\times S^n$:
\begin{itemize}
\item When $n$ is even, the diagonal $D$ has self-intersection number 2 and the anti-diagonal 
$A$ has self-intersection number $-2$.
\item When $n$ is odd, $D$ and $A$ each have self-intersection number 0.
\end{itemize}

\medskip

\noindent\bf Lemma F\it.\ \ Regardless of the parity of $n$, the anti-diagonal $A' = p(A)$ in \linebreak
$(S^n \dst  S^n) - D'$ has self-intersection number \rm 1 mod 2 .

\medskip

When $n = 1$ this is easy to see visually, since $S^1\dst  S^1$ is a M\"obius band, while $A'$ 
is the circle running along the middle of the band.

To prove the lemma in general, we will describe a concrete perturbation of $A'$ in 
$(S^n \dst  S^n) - D'$ which meets $A'$ transversally in just one point.

To that end, let $f\colon S^n\to S^n$ be a diffeomorphism with fixed points at the north and 
south poles, but otherwise moving each point of $S^n$ slightly southwards along its circle of 
longitude We want to choose $f$ to satisfy the following two conditions:
\begin{enumerate}
\item[(1)] The differential of $f$ at the north pole is expansive, and at the south pole contractive.
\item[(2)] The behavior of $f$ is related to the antipodal map as in the following diagram, which 
shows a typical great circle of longitude.
\end{enumerate}

\begin{figure}[h!]
\center{\includegraphics[height=155pt]{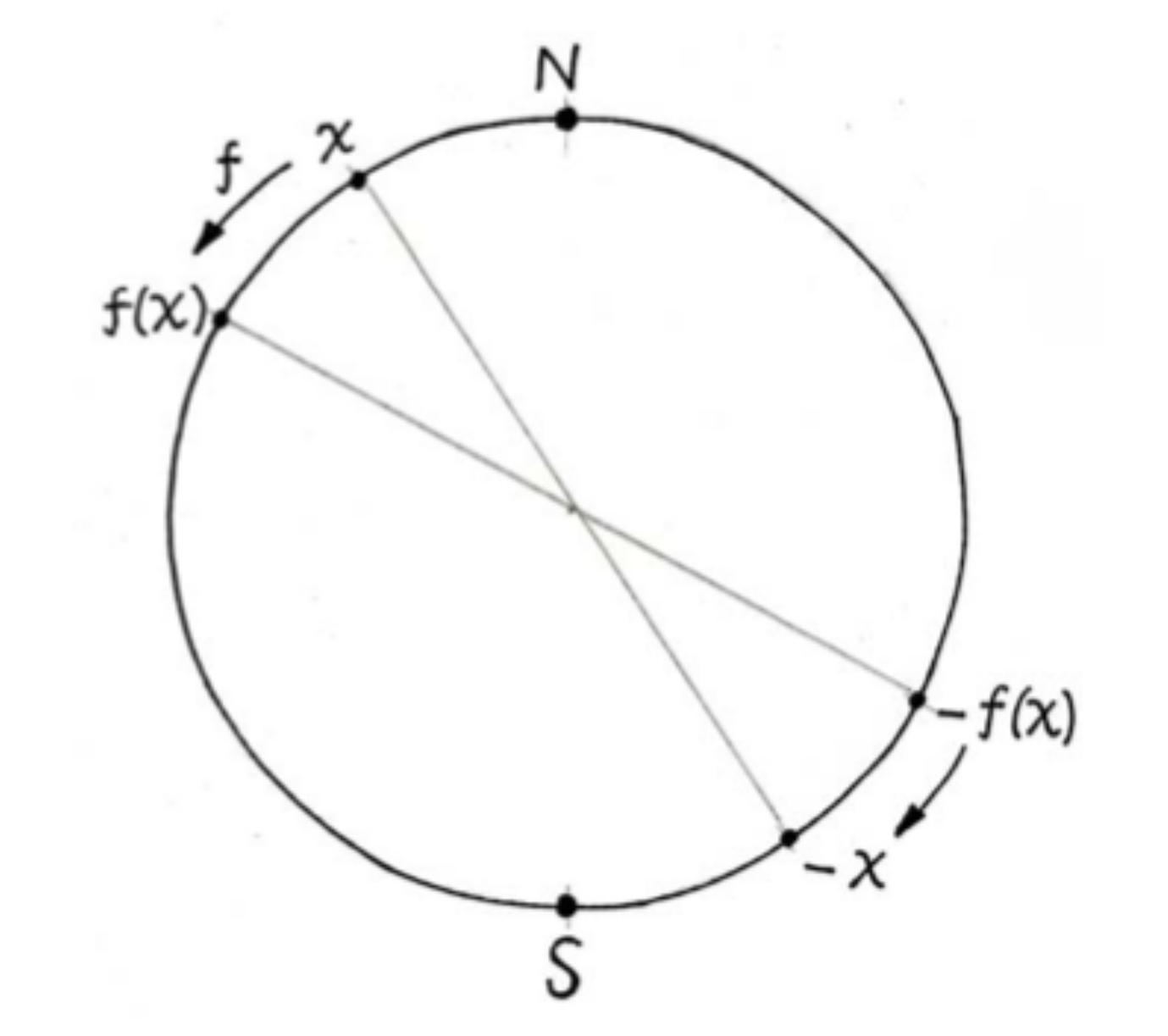}
\caption{\large{\textit{\textbf{\boldmath Required behavior of $f$}}}
}}
\end{figure}

We are requiring that $f(-f(x)) = -x$, or in other words, that the map $-f$ be an involution.

We intend that $f\colon S^n\to S^n$ should be the same on every great semi-circle of longitude, and
construct such a map as follows.

First we redraw the above circle of longitude on $S^n$, focus on its ``left half'', and parametrize 
this from 0 in the north to 1 in the south.

\begin{figure}[h!]
\center{\includegraphics[height=155pt]{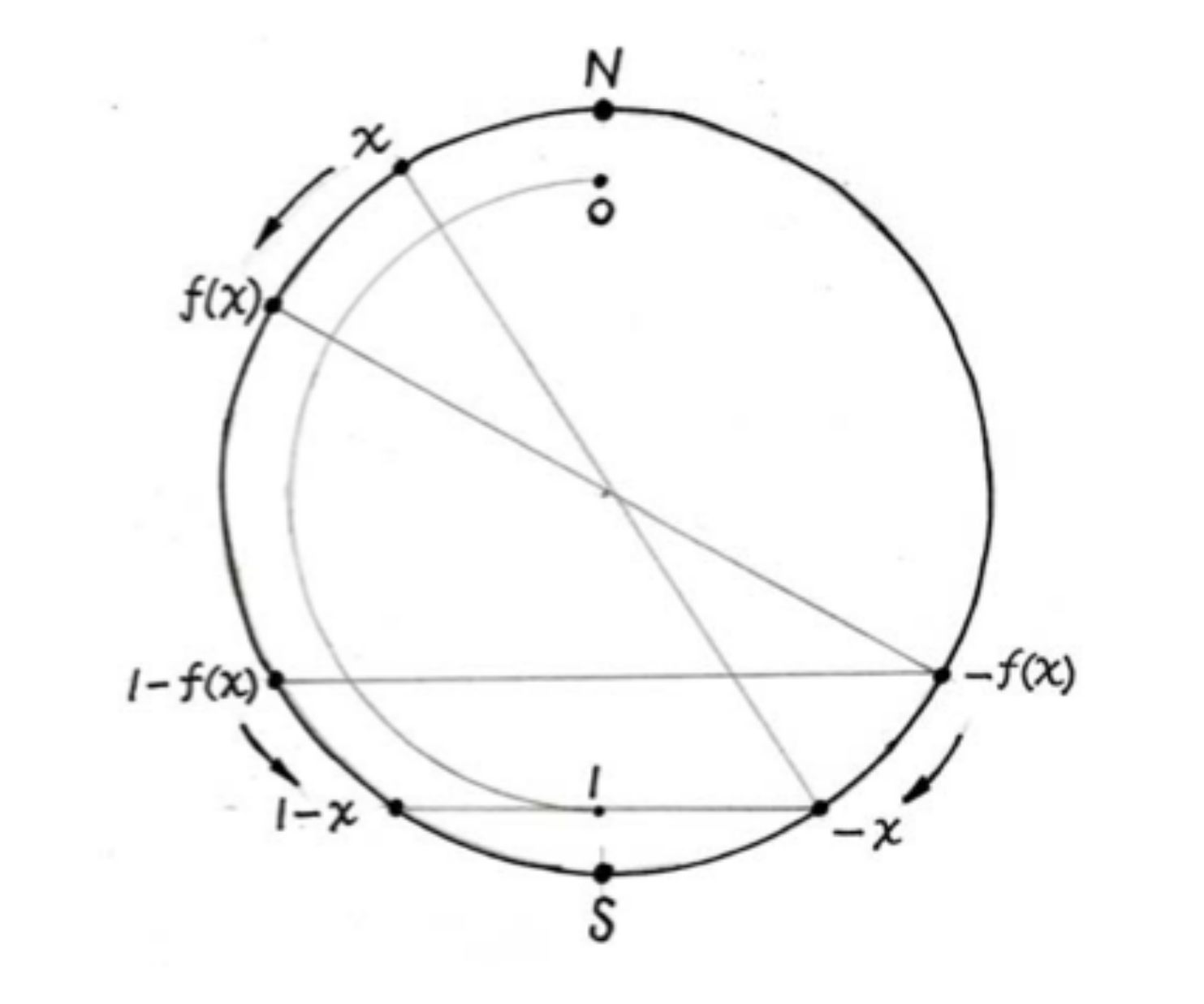}
\caption{\large{\textit{\textbf{\boldmath Behavior of reparametrized $f$}}}
}}
\end{figure}

We see in the above diagram that the point $-f(x)$ corresponds to the point $1- f(x)$ in this 
parametrization, and likewise the point $-x$ corresponds to $1- x$.

So, focusing on the left semi-circle, and thinking of $f$ as a map from $[0, 1]$ to itself, we are 
requiring in condition (2) above that
\begin{enumerate}
\item[(2')]  $f(1-f(x)) = 1-x$.
\end{enumerate}
To construct such a function $f$, we are guided by the following diagram, in which we 
show the graph of $f$ inside the square $[0, 1] \times [0, 1]$.

\begin{figure}[h!]
\center{\includegraphics[height=200pt]{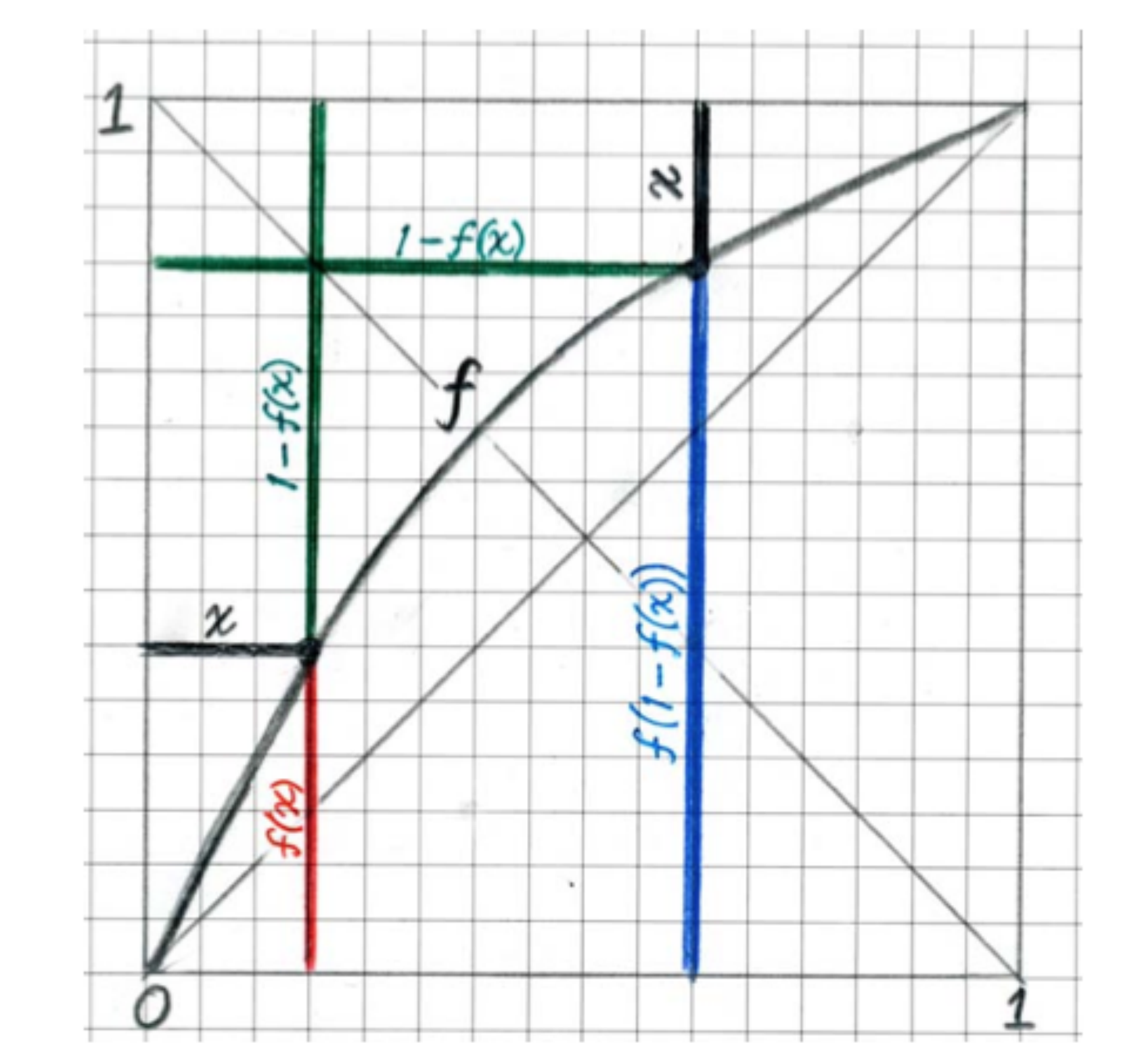}
\caption{\large{\textit{\textbf{\boldmath Guide for constructing $f$}}}
}}
\end{figure}

We insist that the graph of $f$ should be invariant under reflection about the diagonal of slope $-1$, 
and will show that this guarantees condition (2').

To that end, start with the point $(x, f(x))$ on the graph. In the diagram, we show $x$ as a 
horizontal segment in black, $f(x)$ as a vertical segment in red, and then $1- f(x)$ as a vertical 
segment in green above it.

Now reflect in the diagonal of slope $-1$. Then $1-f(x)$ reappears as a horizontal segment in 
green, and due to the symmetry, its right hand end point is still on the graph of $f$.

Hence $f(1 - f(x))$ appears as a vertical segment in blue. 

Then $1 - f(1 - f(x))$ appears directly above it as a vertical segment in black
which, thanks again to the reflective symmetry, has length $x$. That is, $1-f(1-f(x)) = x$,
which is just a transposed version of condition (2').

To take care of condition 1, we simply make $f'(0) > 1$. Then by the reflective symmetry we have 
$f'(1) = 1 / f'(0) < 1$. For example, in Figure 12, we have $f'(0) = 2$ and $f'(1) = 1/2$.

We now define $f\colon S^n\to S^n$ by copying the map $f$ on every great semi-circle of 
longitude from the north pole $N$ to the south pole $S$.

\vfill
\eject

Differentiability of this $f$ at $N$ and $S$ is easily guaranteed by making the map 
$f\colon [0, 1] \to [0, 1]$ linear near $0$, and hence by reflection, also near $1$, as shown in Figure 12.

With such a map $f\colon S^n\to S^n$ in hand, we complete the proof of Lemma F as follows.

The anti-diagonal $n$-sphere $A = \{(x, -x)\}$ in $S^n\times S^n$ projects down by $p$ to the 
anti-diagonal real projective $n$-space $A' = \{[x, -x]\}$ in $S^n \dst  S^n$.

Up in $S^n \times S^n$, consider the smooth $n$-sphere $A_f = \{(x, -f(x)\}$. If $f$ is close to the 
identity, then $A_f$ is a slight perturbation of $A$, which by condition (1) above meets $A$ transversally 
in the two points $(N, S)$ and $(S, N)$.

Furthermore, by condition (2), for each point $(x, -f(x))$ of $A_f$, the point $(-f (x), x)$ is also in $A_f$,
and so $A_f$ is invariant under the involution $\sigma$.

Therefore the map $p\colon A_f\to p(A_f) = A'_f$ is a double covering, whose image is a 
slight perturbation of $A'$, and which meets it transversally in the single point \linebreak
$[N, S] = [S, N]$. 
By construction of $f$, this perturbation of $A'$ takes place entirely in $(S^n \dst  S^n) - D'$.

This completes the proof of Lemma F.

\medskip

\noindent\bf Proof of the Key Lemma.\rm

To prove the Key Lemma as rephrased earlier, we start with two maps $f_1$ and $f_2\colon S^n\to S^n$,
both homotopic to the identity, and must find points $u$ and $v$ in $S^n$ such that $f_1(u)$ and 
$f_1(v)$ are antipodal, and at the same time $f_2(u)$ and $f_2(v)$ are also antipodal.

To that end, define $F_1$ and $F_2\colon S^n \dst  S^n \to S^n \dst  S^n$ by 
$$F_1([u, v]) = [f_1(u), f_1(v)]\quad\mbox{\rm and}\quad	F_2([u, v]) = [f_2(u), f_2(v)].$$
We will show that
$$F_1^{-1}(A') \cap F_2^{-1}(A')\ne\emptyset.$$

Suppose, to the contrary, that the inverse images $F_1^{-1}(A')$ and $F_2^{-1}(A')$ are disjoint. 
Choose an open tubular neighborhood $N'$ of $A'$ in $S^n \dst  S^n - D'$ such that
$$F_1^{-1}(N')\cap F_2^{-1}(N') = \emptyset.$$
In the following argument, all homology and cohomology will be understood to have $\ints_2$	
coefficients. The symbol  $H^{2n}_c$ denotes cohomology with compact supports.

\vfill
\eject

Using the terminology and results of Eilenberg-Steenrod XI.6 , the pair\linebreak
$(S^n \dst  S^n, D')$ is a 
relative $2n$-manifold. The fact that $S^n \dst  S^n	- D'$ is connected implies that
$$H^{2n}(S^n \dst  S^n,D') \cong \ints_2,\quad H^{2n}_c(S^n \dst  S^n -D') \cong \ints_2,$$ 
and the inclusion $(S^n \dst  S^n - D', \emptyset) \subset (S^n \dst  S^n, D')$ induces an isomorphism
$$H^{2n}_c(S^n \dst S^n -D') \to H^{2n}(S^n \dst S^n,D').$$

Consider the compact smooth submanifold $A' \subset S^n \dst  S^n - D'$. Let $\omega$ be an 
$n$-dimensional cochain representing the Poincar\'e dual of $A'$ in the ring 
$H^{2n}_c(S^n \dst  S^n - D')$, that is, $\omega$ represents the Thom class of the normal bundle 
of $A'$. We may assume that $\omega$ is supported in the open tubular neighborhood $N'$ of $A'$.

We saw in the previous section that $A'$ has nonzero self-intersection number mod 2 in 
$S^n \dst  S^n - D'$. It follows from Poincar\'e duality that the cup product $[\omega] \cup [\omega]$ 
is nonzero in $H^{2n}_c(S^n \dst  S^n - D')$. Using the inclusion map, we see that 
$[\omega] \cup [\omega]$ is also nonzero in $H^{2n}(S^n \dst  S^n, D')$.

Since $f_1$ and $f_2\colon S^n\to S^n$ are both homotopic to the identity, it follows that \linebreak
$F_1$ and $F_2\colon S^n \dst  S^n \to S^n \dst  S^n$ are also both homotopic to the identity through 
maps which always take the singular locus $D'$ to itself (though not keeping it pointwise fixed). 
Therefore
$$F_1^\dast[\omega] \cup F_2^\dast[\omega] = [\omega] \cup [\omega] = 1$$ 
in $H^{2n}(S^n \dst S^n,D') \cong \ints_2$.

But since the support of $\omega$ is contained in $N'$, and we have seen that \linebreak
$F_1^{-1}(N') \cap F_2^{-1}(N') = \emptyset$, it follows that the cochains $F_1^\dast\omega$
and $F_2^\dast\omega$ have disjoint supports. Therefore
$$F_1^\dast[\omega] \cup F_2^\dast[\omega] = 0.$$ 
This contradiction shows that
$$F_1^{-1}(A') \cap F_2^{-1}(A')\ne\emptyset.$$ 
Now let $[u, v]$ be a point in this intersection.

The fact that $F_1([u, v]) = [f_1(u), f_1(v)]$ lies in $A'$ tells us that $f_1(u)$ and $f_1(v)$ are antipodal.

The fact that $F_2([u, v]) = [f_2(u), f_2(v)]$ lies in $A'$ tells us that $f_2(u)$ and $f_2(v)$ are antipodal.

This is exactly the claim of the Key Lemma, and so completes its proof, and with it the proof of Theorem B.

\vfill
\eject

\centerline{\bf PART IV. PROOF OF THEOREM C} 

\medskip

\noindent\bf Equivalence of Hopf vector fields.\rm

We will make use of the following in the proof of Theorem C.

An \it orthogonal complex structure \rm $J$ on Euclidean space $\reals^{2n+2}$ is an element of 
$SO(2n+2)$ such that $J^2 = -\,\mbox{\rm Identity}$. Decomposing $\reals^{2n+2}$ 
into an orthogonal direct sum of 2-planes invariant under $J$, we note that
\begin{enumerate}
\item[(i)] $x\cdt J(x) = 0$ for all $x\in\reals^{2n+2}$, and 
\item[(ii)] Any two orthogonal complex structures on $\reals^{2n+2}$
are conjugate in $O(2n+2)$.
\end{enumerate}
A Hopf vector field on $S^{2n+1}$ is the same thing as the restriction to $S^{2n+1}$ 
of an orthogonal complex structure on $\reals^{2n+2}$.

Suppose that $J$ and $J'$ are any two orthogonal complex structures on $\reals^{2n+2}$, and that $g$
is an element of $O(2n+2)$, thanks to (ii) above, such that $g J = J' g$. Then $(g, g)$ is an 
isometry of the unit tangent bundle $US^{2n+1}$ taking the graph $V$ of the restriction of $J$ to 
$S^{2n+1}$ to the graph $V'$ of
the corresponding restriction of $J'$, since 
$$(g, g)(x, J(x)) = (g(x) , gJ(x)) = (g(x) , J'g(x)) = (y, J'(y).$$ 
In other words, any two Hopf cross sections $V(S^{2n+1})$ and $V'(S^{2n+1})$ of 
$US^{2n+1}$ can be taken to one another by an isometry of this unit tangent bundle.

\medskip

\noindent\bf Proof of Theorem C.\rm

Let $v$ be a Hopf vector field on $S^{2n+1}$, obtained as the restriction of the orthogonal 
complex structure $J$ on $R^{2n+2}$. Let $V$ be the corresponding cross-section of 
$US^{2n+1}$, and denote by $i\colon V(S^{2n+1})\to US^{2n+1}$ the inclusion map. The claim of Theorem C is that if the map $f\colon V(S^{2n+1}) \to US^{2n+1}$ is 
homotopic to the inclusion and has Lipschitz constant $\le 1$, then its Lipschitz constant 
equals 1 and $f$ is isometric to the inclusion.

The composite inclusion $V(S^{2n+1}) \subset US^{2n+1}\subset S^{2n+1} \times S^{2n+1}$ 
is isometric to the inclusion of the diagonal $\Delta S^{2n+1} \subset S^{2n+1}\times S^{2n+1}$, 
since the restriction of $J$ to $S^{2n+1}$ is an isometry.

The composite map $f\colon V(S^{2n+1}) \to US^{2n+1}\subset S^{2n+1}\times S^{2n+1}$ is 
homotopic to the inclusion of $V(S^{2n+1})$ into $S^{2n+1} \times S^{2n+1}$, and still has 
Lipschitz constant $\le 1$ there. So by Theorem B, the composite map $f$ must have Lipschitz constant equal to 1 and be 
isometric to the inclusion of $V(S^{2n+1})$ into $S^{2n+1} \times S^{2n+1}$.

Thus the image under $f$ of $V(S^{2n+1})$ in $S^{2n+1} \times S^{2n+1}$ must be the 
graph $V'$ of an orientation-preserving isometry $J'\colon S^{2n+1}\to S^{2n+1}$ such 
that $x \cdt J'(x) = 0$ for all $x$ in $S^{2n+1}$. In other words, $f$ takes the Hopf 
cross-section $V(S^{2n+1})$ to another Hopf cross-section $V'(S^{2n+1})$ in $US^{2n+1}$.

Since both $V(S^{2n+1})$ and $V'(S^{2n+1})$ are round $2n+1$-spheres of radius $\sqrt{2}$,
and $f$ has Lipschitz constant $\le 1$, the map $f\colon V(S^{2n+1})\to V'(S^{2n+1})$ must be 
an isometry.

We saw in the previous section that there is an isometry of $US^{2n+1}$ to itself which takes 
$V(S^{2n+1})$ to $V'(S^{2n+1})$. It follows that $f\colon V(S^{2n+1}) \to US^{2n+1}$ is isometric 
to the inclusion $i\colon V(S^{2n+1}) \to US^{2n+1}$, completing the proof of Theorem C.

So we see that Theorem C is a direct consequence of Theorem B.

\medskip

\noindent\bf Comment on Theorem C.\rm

There are really \it two \rm distinct Riemannian metrics on the unit tangent bundle of a sphere.

The first, which we have been using, views 
$$US^n = \{(x,y)\,:\,x,y\in S^n\ \mbox{\rm and}\ x\cdt y=0\} \subset S^n \times S^n,$$
takes the usual product metric on $S^n \times S^n$, and then gives to $US^n$ the 
Riemannian metric induced by this inclusion.

The second, due to Sasaki [1958], applies to the tangent bundle $TM$ of any Riemannian 
manifold $M$. If $(x(t), v(t))$ is a curve in $TM$, then the length of the tangent vector 
to this curve is taken to be
$$(|x'(t)|^2 + |v'(t)|^2)^{1/2},$$
where $x'(t)$ is the tangent vector to the curve $x(t)$ in $M$, where $v'(t)$ is the covariant 
derivative of the vector field $v(t)$ along the curve $x(t)$, and the norms of these vectors 
are measured in the given Riemannian metric on $M$.

As a result, if $v(t)$ is a parallel vector field along the curve $x(t)$ in $M$, meaning that the 
covariant derivative $v'(t) = 0$, then the length of the tangent vector to the curve $(x(t), v(t))$ 
in $TM$ is simply the length $|x'(t)|$ of the tangent vector to the curve $x(t)$ in $M$.

For example, if $M$ is the unit circle $S^1$ in $\reals^2$, and if 
$$x(t) = (\cos t,\sin t)\quad\mbox{\rm and}\quad v(t) = (-\sin t,\cos t),$$
then not only is $x(t)$ a unit speed curve in $M$, but also $(x(t), v(t))$ is a unit speed curve in $TM$,
since $v(t)$ is parallel along $x(t)$. In other words, the fact that $v(t)$ is, to the naked eye, spinning around just as fast as $x(t)$, 
is forgiven, and the length of the loop $(x(t), v(t))$ in $TM$ is just $2\pi$, as opposed to $2\pi\sqrt{2}$.

The same thing happens with the Sasaki metric on the unit tangent bundle $US^n$ of $S^n$. 
If $x(t)$ travels at unit speed once around a great circle in $S^n$, and if $v(t) = x'(t)$ is its 
velocity vector, then the curve $(x(t), v(t))$ in $US^n$ also travels at unit speed in $US^n$, and so 
has length $2\pi$. By contrast, in the ``product metric'' on $US^n$ inherited from its natural 
inclusion in $S^n\times S^n$, this loop has length $2\pi\sqrt{2}$.

This is the only difference between the two competing metrics on $US^n$: you pass from the 
product metric to the Sasaki metric by reducing lengths by a factor of $\sqrt{2}$ in 
the direction of the above ``geodesic flow'', while preserving lengths in the orthogonal direction.

It would be sensible to check the validity of Theorem C using the Sasaki metric on the 
unit tangent bundle $US^{2n+1}$, but we have not done this yet.

\vfill
\eject

\centerline{\bf PART V. PROOF OF THEOREM D.}

\medskip

\noindent\bf Set up.\rm

The Stiefel manifold $V_2\reals^4$ is the set of orthonormal 2-frames in 4-space, 
$$V_2\reals^4 = \{(x,y)\,:\, x\in S^3,\ y\in S^3,\ x\cdt y=0\},$$
and we give it the Riemannian metric inherited from its inclusion in $S^3\times S^3$.

The Grassmann manifold $G_2\reals^4$ is the set of oriented 2-planes through the origin in 4-space. 
We identify it with the set of unit decomposable 2-vectors in the exterior product $\Lambda^2\reals^4$, 
a 6-dimensional Euclidean space, and give it the resulting Riemannian metric.

\medskip

\noindent\bf Fact \rm (Gluck and Warner [1983])\it.\ \  A $2$-vector in $4$-space is decomposable if 
and only if it has equal length projections into the $+1$ and $-1$  eigenspaces $E^3_+$	
and $E^3_-$ of the Hodge star operator on $\Lambda^2\reals^4$.\rm

\medskip

From this fact, it follows that 
$$G_2\reals^4 = S^2_+(1/\sqrt{2}) \times S^2_-(1/\sqrt{2}) \subset E^3_+\oplus E^3_- 
= \Lambda^2\reals^4,$$
the product of the 2-spheres of radius $1/\sqrt{2}$ in $E^3_+$ and $E^3_-$.

The projection map $p\colon V_2\reals^4\to G_2\reals^4$ takes 
$(x, y) \to (x\wedge y)/\Vert x \wedge y\Vert.$ It is a Riemannian submersion, and thus has 
Lipschitz constant 1 .

Theorem D asserts that any map homotopic to the Stiefel projection has Lipschitz constant 
$\ge 1$, with equality if and only if the map is isometric to this projection.

In other words, the Stiefel projection is, up to isometries of domain and range, the unique Lipschitz 
constant minimizer in its homotopy class.

To prove this theorem, we will observe within the Stiefel projection $V_2\reals^4\to G_2\reals^4$ 
two families of Hopf projections $S^3\to S^2$, whose Lipschitz minimality, unique up to 
isometries of domain and range, was established in Theorem A. They provide the framework for 
the proof.

\vfill
\eject

\noindent\bf An alternative view of the Stiefel projection.\rm

The Stiefel manifold $V_2\reals^4$ is the same as the unit tangent bundle $US^3$ of the 3-sphere. 
This bundle is trivial topologically (though not metrically), and has two common sense 
trivializations, $US^3\to S^3\times S^2$, given by
$$(x,y) \to (x,yx^{-1})\quad\mbox{\rm and}\quad (x,y) \to (x,x^{-1}y),$$ 
using multiplication of unit quaternions, and thinking of $S^2$ as the space
of purely imaginary unit quaternions. 

Packaging these two trivializations together yields a map
$$V_2\reals^4 = US^3\to S^2\times S^2 \quad \mbox{\rm by}\quad (x,y)\to (yx^{-1},x^{-1}y),$$ 
which a simple computation shows to be just a copy of the Stiefel projection\linebreak
$p\colon V_2\reals^4	\to G_2\reals^4$, scaled up by the linear factor $\sqrt{2}$.

This version of the Stiefel projection $p$ has Lipschitz constant $\sqrt{2}$, and we will use 
it in what follows.

\medskip

\noindent\bf Copies of the complex Grassmannian $G_1\cplxs^2$ inside $G_2\reals^4$.\rm

It is easy to see that on $\reals^4$, all orthogonal complex structures are given by left or 
right multiplication by a purely imaginary unit quaternion. To be specific, let us use left multiplication 
by the purely imaginary unit quaternion $u$ to regard $\reals^4$ as $\cplxs^2$.

The corresponding complex Grassmannian $G_1\cplxs^2$ consists of all complex lines in 
$\cplxs^2$ through the origin. To real eyes, each such complex line is a 2-plane through 
the origin, with a natural orientation given by the ordered basis $x,\, ux$ for any unit vector $x$ therein.

Using the alternative version $p(x, y) = (yx^{-1},x^{-1}y)$ of the Stiefel projection, we have
$p(x, ux) = (u, x^{-1}ux)$, which tells us that
$$G_1\cplxs^2 = u\times S^2 \subset S^2\times S^2 = G_2\reals^4,$$ 
a ``vertical'' 2-sphere in $S^2\times S^2$.
The inverse image $p^{-1}(G_1\cplxs^2)$ up in $V_2\reals^4$ is the subset 
$${}_uV = \{(x,ux)\,:\, x\in S^3\}.$$
It is a round, totally geodesic 3-sphere of radius $\sqrt{2}$, and is the graph of the 
corresponding Hopf vector field.

The restriction $p\colon {}_uV \to G_1\cplxs^2$ of the Stiefel projection is just a copy of the 
Hopf projection, scaled up by a factor $\sqrt{2}$.

Varying the choice of purely imaginary unit quaternion $u$ gives us all possible 
``vertical'' 2-spheres $u \times S^2$ as the corresponding $G_1\cplxs^2$ inside our 
$S^2 \times S^2$ picture of $G_2\reals^4$. 
Each one serves as the base space of a Hopf projection, as above.

Similarly, fixing an orthogonal complex structure on $\reals^4$ via right multiplication by the 
purely imaginary unit quaternion $v$, we get the corresponding
$$G_1\cplxs^2 = S^2\times v \subset S^2\times S^2 = G_2\reals^4,$$ 
a ``horizontal'' 2-sphere in $S^2\times S^2$.

The inverse image $p^{-1}(G_1\cplxs^2)$ up in $V_2\reals^4$ is the subset 
$$V_v = \{(x,xv)\,:\, x\in S^3\},$$
and again, the restriction $p\colon V_v\to G_1\cplxs^2$ of the Stiefel projection is a 
scaled-up copy of the Hopf projection.

Varying the choice of $v$ gives us all possible ``horizontal'' 2-spheres $S^2\times v$ as 
the corresponding $G_1\cplxs^2$ inside our $S^2\times S^2$ picture of $G_2\reals^4$, again 
with each serving as the base space of a Hopf projection.

\begin{figure}[h!]
\center{\includegraphics[height=220pt]{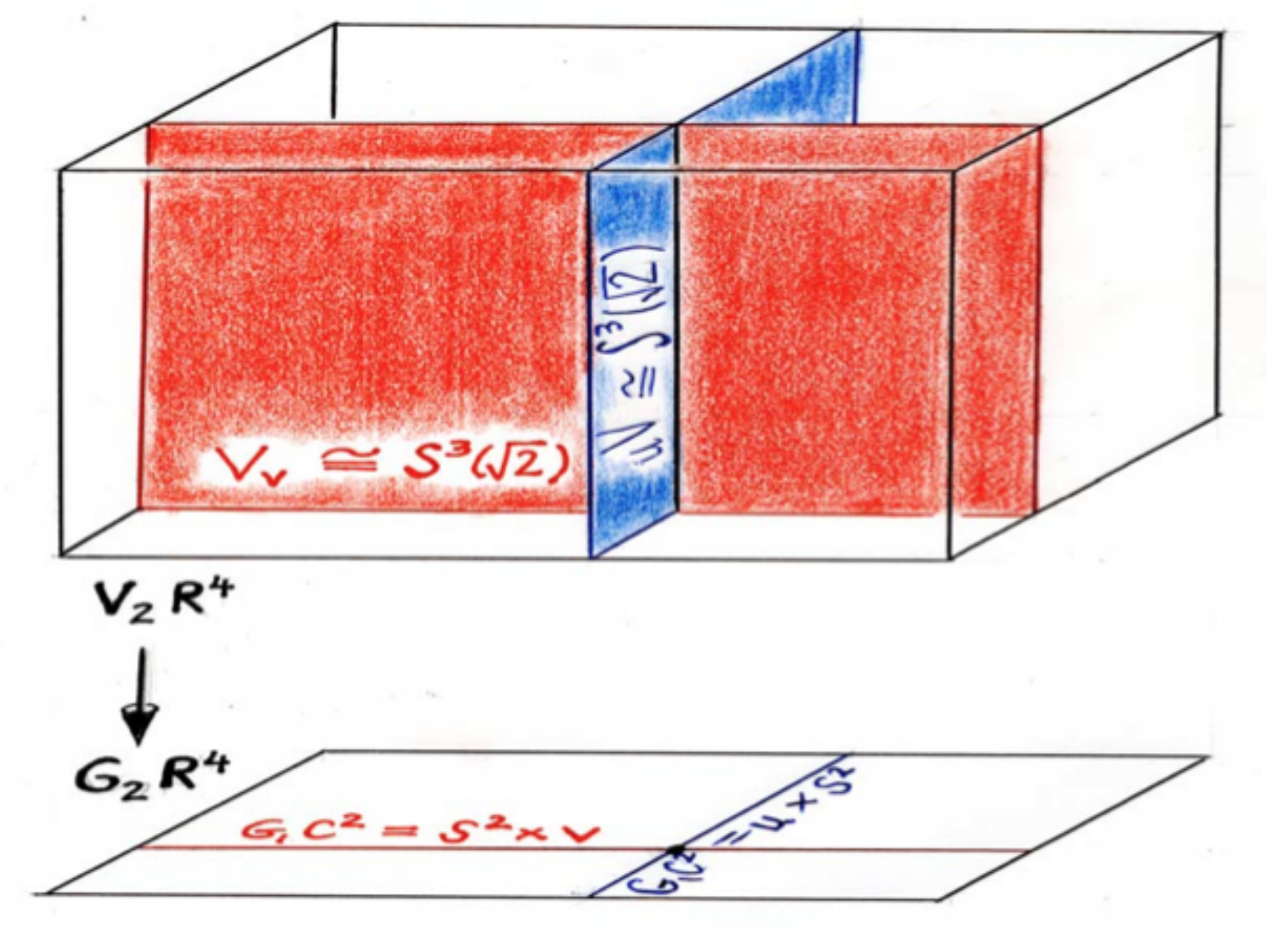}
\caption{\normalsize{\textit{\textbf{\boldmath Two copies of $G_1\cplxs^2$ inside 
$G_2\reals^4$, and their inverse images in 
$V_2\reals^4$}}}}
}
\end{figure}

\vfill
\eject

\noindent\bf Proof of Theorem D.\rm

We start with the Stiefel projection $p\colon V_2\reals^4\to G_2\reals^4= S^2\times S^2$, 
a scaled-up Riemannian submersion with Lipschitz constant $\sqrt{2}$.

Then we consider another map $f\colon V_2\reals^4\to G_2\reals^4$ which is homotopic to 
$\,p\,$ with Lipschitz constant $\le \sqrt{2}$, and must show that $f$ has Lipschitz constant equal 
to $\sqrt{2}$, and agrees with $\,p\,$ up to composition with isometries of domain and range.

The isometries of the domain $V_2\reals^4$ are known to preserve the Stiefel fibres 
(Gluck and Ziller [1986]), so somewhere in the argument we will have to show that these are also 
the fibres of $f$, and we take this as a hint.

Fix a purely imaginary unit quaternion $v$, and use right multiplication by $v$
to impose an orthogonal complex structure on $\reals^4$, so that we can regard it as $\cplxs^2$.

Then consider again the restriction $p\colon V_v\to G_1\cplxs^2 \subset G_2\reals^4$ 
of the Stiefel projection.

Let $p_1$ and $p_2$ denote the compositions of $p$ with the projections of 
$G_2\reals^4= S^2\times S^2$ to its first and second factors.

Likewise, let $f_1$ and $f_2$ denote the corresponding compositions of $f$ with these two 
projections.

Now compare the restrictions 
$$p_1\ \mbox{\rm and}\  f_1\colon V_v \cong S^3(\sqrt{2})\to S^2.$$
Since $p$ and $f$ are homotopic, so are the above restrictions of $p_1$ and $f_1$. 
The restricted $p_1$ is a Hopf projection, scaled up by $\sqrt{2}$, with Lipschitz 
constant $\sqrt{2}$, while the restricted $f_1$ has Lipschitz constant $\le\sqrt{2}$.

By Theorem A, the restricted $f_1$ must have Lipschitz constant equal to $\sqrt{2}$, and agree 
with the restricted $p_1$ up to isometries of domain and range.

In other words, $f_1\colon V_v	\cong S^3(\sqrt{2})\to S^2$ is a Hopf projection, 
scaled up by $\sqrt{2}$.

\medskip

\noindent\bf Claim\it.\ \  The image $f(V_v)$ is a horizontal 2-sphere 
$S^2\times v'$ in $S^2\times S^2$.\rm

\noindent\it Proof of claim\rm.\ \ We already know that $f_1\colon V_v\cong S^3(\sqrt{2})\to S^2$ is 
a Hopf projection, whose fibres are great circles of radius $\sqrt{2}$.

In particular, the map $f_1$ is smooth.

In the tangent space to $V_v$ at the point $(x, xv)$, let $F_x$ denote the tangent line to the 
great circle fibre of $f_1$, and $G_x$ the orthogonal 2-plane.

The differential of $f_1$ at this point takes $G_x$ conformally to the tangent 2-plane to $S^2$ at 
the image point, stretching lengths by $\sqrt{2}$.

In particular, any smooth curve in $V_v$ everywhere tangent to $G_x$ is taken by $f_1$ conformally 
to a curve in $S^2$, stretching lengths by $\sqrt{2}$.

Since this is the maximum stretch allowed the map $f\colon V_v\to G_2\reals^4 = S^2\times S^2$, 
that same curve in $V_v$ must be taken by $f$ to a horizontal curve in $S^2\times S^2$.

But any two Hopf fibres can be connected by a smooth curve in $S^3$ which is 
everywhere orthogonal to the fibres its passes through, and hence $f$ must take all of $V_v$ 
to a single horizontal 2-sphere $S^2\times v'$, verifying the claim.

And since $f_1\colon V_v\to S^2$ is a Hopf fibration, so also is $f\colon V_v \to S^2\times v'$, 
where it is perfectly possible that $v'\ne v$.

We have thus gained some control over the nature of $f$: on each $V_v$ the map $f$ must be a 
Hopf projection, with image a horizontal 2-sphere in $S^2 \times S^2$.

The fibres of $f$ must therefore be great circles in $V_v\cong S^3(\sqrt{2})$, but we don't yet know that they coincide with the Stiefel fibres.

Now repeat all of the above with orthogonal complex structures on $\reals^4$ given by 
\it left \rm multiplication by a purely imaginary unit quaternion $u$, and learn that the image 
$f({}_uV)$ is a \it vertical \rm 2-sphere $u'\times S^2$ in $S^2\times S^2$, where again it is possible 
that $u'\ne u$.

The two images $u'\times S^2$ and $S^2\times v'$ intersect in the single point $(u', v')$.

Since the inverse image of an intersection is the intersection of the inverse images, the portion 
of $f^{-1}(u', v')$ within ${}_uV \cup V_v$ must be ${}_uV\cap V_v$, which is a Stiefel fibre.

It follows in this way that the map $f$ must collapse each Stiefel fibre in $V_2\reals^4$ to a single 
point in $G_2\reals^4$, although we do not yet know that distinct Stiefel fibres are sent to 
distinct points by $f$.

Since $f \colon V_2\reals^4\to G_2\reals^4$ collapses each Stiefel fibre to a point, it induces
a map
$$\bah{f}\colon G_2\reals^4 =S^2\times S^2\to G_2\reals^4 =S^2\times S^2$$
with Lipschitz constant $\le 1$.

From the above discussion, we see that $\bah{f}$ takes the horizontal 2-sphere $S^2\times v$ to 
the horizontal 2-sphere $S^2\times v'$, and the vertical 2-sphere $u\times S^2$	 to the vertical 
2-sphere $u'\times S^2$, in each case with Lipschitz constant $\le 1$.

Since by hypothesis the map $f$ is homotopic to the Stiefel projection $p$, these 
horizontal-to-horizontal and vertical-to-vertical maps of 2-spheres are all homotopic to the identity.

But a map from $S^2$ to $S^2$ which is homotopic to the identity and has Lipschitz constant $\le 1$
must be an orientation-preserving isometry.

It follows that
$$\bah{f}(u, v) = (u', v') = (g(u), h(v)),$$ 
where $g$ and $h$ are orientation-preserving isometries of $S^2$.

Hence $f$ differs from the Stiefel projection $p$ by an isometry $(g, h)$ of their common range 
$G_2\reals^4$, completing the proof of Theorem D.

\vfill
\eject

\centerline{\bf REFERENCES\rm}

\addtolength{\baselineskip}{-2pt}

\begin{description}
\item[1931] Heinz Hopf, {\it Uber die Abbildungen der dreidimensionalen Sph\"are auf die 
K\"ugelflache}, Math.\ Ann.\ {\bf 104}, 637--665.

\item[1935] Heinz Hopf, {\it Uber die Abbildungen von Sph\"aren auf Sph\"aren niedrigerer 
Dimension}, Fund.\ Math.\ {\bf 25}, 427--440.

\item[1947]  J.H.C.~Whitehead, {\it An expression of Hopf's invariant as an integral},
Proc.\ Natl.\ Acad.\ Sci.\ USA {\bf 33}(5), 117--123. 

\item[1958] Shigeo Sasaki, {\it On the differential geometry of tangent bundles of 
Riemannian manifolds}, Tohoku Math. J.\ (2) {\bf 10}, 338--354.

\item[1960] Marcel Berger, {\it Les var\'et\'es riemanniennes $1/4$-pinc\'ees}, Ann.\ Scu.\ Norm.\ Sup.\ Pisa {\bf 14}, 161--170.

\item[1961] Y.-C.\ Wong, {\it Isoclinic $n$-planes in Euclidean $2n$-space, Clifford parallels in elliptic
$2n-1$ space, and the Hurwitz matrix equations}, Memoirs Amer.\ Math.\ Soc.\ {\bf 41}, 1--112.

\item[1963] William Browder, {\it Higher torsion in $H$-spaces}, Trans.\ Amer.\ Math.\  Soc.\  {\bf 108}, 353--375.

\item[1963] Joseph Wolf, {\it Geodesic spheres in Grassmann manifolds}, Illinois J.\ Math.\ {\bf7}, 
425--446.

\item[1963] Joseph Wolf, {\it Elliptic spaces in Grassmann manifolds}, Illinois J.\  Math.\  {\bf 7}, 
447--462.

\item[1966] Barrett O'Neill, {\it The fundamental equations of a submersion}, Mich.\ Math.\ J.\ {\bf 16}, 
459--469.

\item[1975] Richard H.\ Escobales, {\it Riemannian submersions with totally geodesic fibers}, 
J.\  Diff.\ Geom.\ {\bf 10}, 253--276.

\item[1975] R.T.\ Smith, {\it Harmonic mappings of spheres}, Amer.\ J.\ Math.\ {\bf 97}, 364--385.

\item[1977] Odette Nouhaud, {\it Applications harmoniques d'une vari\'et\'e riemanienne dans son 
fibr\'e tangent. Generalisation}, C.R.\ Acad.\ Sci.\ Paris, S\'er.\ A-B {\bf 284}, No.\ 14, A815--A818.

\item[1978] James Eells and Luc Lemaire, {\it A report on harmonic maps}, Bull.\ London 
Math.\ Soc.\ {\bf 10}(1), 1--68.

\item[1979] Toru Ishihara, {\it Harmonic sections of tangent bundles}, J.\ Math.\ Tokushima 
Univ.\ {\bf 13}, 23--27.

\item[1980] Y.L.\ Xin, {\it Some results on stable harmonic maps}, Duke Math.\ J.\  {\bf 47}(3), 609--613.

\item[1982] Raoul Bott and Loring W.\ Tu, {\it Differential forms in algebraic topology}, Graduate Texts 
in Mathematics {\bf 82}, Springer-Verlag, New York.

\item[1983] Herman Gluck and Frank Warner, {\it Great circle fibrations of the three-sphere}, 
Duke Math.\ J.\ {\bf 50}, 107--132.

\item[1985] Detlef Gromoll and Karsten Grove, {\it One dimensional metric foliations in constant 
curvature spaces}, in {\it Differential Geometry and Complex Analysis}, Springer Verlag, 165--168.

\item[1985] Akhil Ranjan, {\it Riemannian submersions of spheres with totally geodesic fibres}, 
Osaka J.\ Math.\ {\bf 22}(2), 243--260.

\item[1985] Shihshu Walter Wei, {\it Classification of stable currents in the product of spheres}, preprint.

\item[1986] Herman Gluck, Frank Warner and Wolfgang Ziller, {\it The geometry of the Hopf 
fibrations}, Enseign.\ Math.\ (2), {\bf 32}(3-4), 173--198.

\item[1986] Herman Gluck and Wolfgang Ziller, {\it On the volume of a unit vector field on the 
three-sphere}, Comment.\ Math.\ Helvetici {\bf 61}(2), 177--192.

\item[1987] Herman Gluck, Frank Warner and Wolfgang Ziller, {\it Fibrations of spheres by 
parallel great spheres and Berger's rigidity theorem}, Ann.\ Global Analysis and Geometry
{\bf 5}(1), 53--82.

\item[1988] Detleff Gromoll and Karsten Grove, {\it The low dimensional metric foliations of 
euclidean spheres}, J.\ Diff.\ Geom.\ {\bf 28}, 143--156.

\item[1988] David L.\ Johnson, {\it Volumes of flows}, Proc.\ Amer.\ Math.\ Soc.\ {\bf 104}(3), 923--931.

\item[1992] Jerzy Konderak, {\it On harmonic vector fields}, Publ.\ Mat.\ {\bf 36}(1), 217--228.

\item[1993] Sharon L. Pedersen, {\it Volumes of vector fields on spheres}, 
Trans.\ Amer.\ Math.\ Soc.\ {\bf 336}, 69--78.

\item[1995] G.\ Wiegemink, {\it Total bending of vector fields on Riemannian manifolds}, 
Math.\ Ann.\ {\bf 303}, 325--344.

\item[1997] Christopher M.\ Wood, {\it On the energy of a unit vector field}, Geom.\ Dedicata 
{\bf 64}(3), 319--330.

\item[1998] Dong-Soong Han and Jin-WhanYim, {\it Unit vector fields on spheres, which are harmonic
maps}, Math.\ Z.\ {\bf 227}(1), 83--92.

\item[2000] Christopher M.\ Wood, {\it The energy of Hopf vector fields},
Manuscripta Math.\ {\bf 101}(1), 71--88.

\item[2000] Fabiano Brito,{\it Total bending of flows with mean curvature correction}, 
Diff.\ Geom.\ Appl.\ {\bf 12}, 157--163.

\item[2000] Fabiano Brito and Pawel Walczak, {\it On the energy of unit vector fields with isolated
singularities}, Ann.\ Math.\ Polon.\ {\bf 73}, 269--274.

\item[2001] Burkhard Wilking, {\it Index parity of closed geodesics and rigidity of Hopf fibrations},
Invent.\ Math.\ {\bf 114}, 281--295.

\item[2001] Olga Gil-Medrano, {\it Relationship between volume and energy of vector fields}, 
Diff.\ Geom.\ Applic.\ {\bf 15}, 137--152.

\item[2002] Olga Gil-Medrano and E.\ Llinares-Fuster, {\it Minimal unit vector fields}, 
Tohoku Math.\ J.\ {\bf 54}, 71--84.

\item[2002] Olga Gil-Medrano, {\it Volume and energy of vector fields on spheres. A survey}, 
in {\it Differential Geometry Valencia} 2001, World Sci.\ Publ., Rivers Edge, NJ, 167--178.

\item[2003] Vincent Borrelli, FabianoBrito and Olga Gil-Medrano, {\it The infimum of the energy of 
unit vector fields on odd-dimensional spheres}, Ann.\ Global Anal.\ Geom.\ {\bf 23}(2), 129--140.

\item[2005] Olga Gil-Medrano, {\it Unit vector fields that are critical points of the volume and of 
the energy: characterization and examples}, in {\it Complex, Contact and Symmetric 
Manifolds}, Birkhauser, Boston, 165--186.

\item[2006] Vincent Borrelli and Olga Gil-Medrano, {\it A critical radius for unit Hopf vector fields 
on spheres}, Math.\ Ann.\ {\bf 334}(4), 731--751.

\end{description}

\noindent University of Pennsylvania\\
Philadelphia, PA  19104

\noindent
\it deturck@math.upenn.edu\rm\\
\it gluck@math.upenn.edu\rm\\
\it peterastorm@gmail.com\rm\\

\end{document}